\documentclass[10pt]{article}
\usepackage{latexsym}
\usepackage{epic}
\usepackage{eepic}
\usepackage{color}
\usepackage{amsmath,multirow,amsfonts}
\usepackage{tikz-cd}
\usetikzlibrary{matrix}
\usepackage{array}
\usepackage[english]{babel}
\usepackage[latin1]{inputenc}
\usepackage{times}
\usepackage[T1]{fontenc}
\usepackage{amssymb}
\usepackage{graphics}
\usepackage{graphicx}
\usepackage{epsfig}
\usepackage{cite}
\usepackage{caption}
\usepackage{tikz}
\usetikzlibrary{arrows}
\tikzstyle{block}=[draw opacity=0.7,line width=1.4cm]
\begin{document}
\newtheorem{def1}{Definition}[section]
\newtheorem{lem}{Lemma}[section]
\newtheorem{exa}{Example}[section]
\newtheorem{thm}{Theorem}[section]
\newtheorem{pro}{Proposition}[section]
\newtheorem{cor}{Corollary}[section]
\newtheorem{rem}{Remark}[section]
\newtheorem{exam}{Example}[section]
\title{$F$-transforms determined by overlap and grouping maps over a complete lattice}
\author{Abha Tripathi\thanks{tripathiabha29@gmail.com}, S.P. Tiwari\thanks{sptiwarimaths@gmail.com}, and {Sutapa Mahato\thanks{sutapaiitdhanbad@gmail.com}}\\
Department of Mathematics \& Computing\\ Indian Institute of Technology (ISM)\\
Dhanbad-826004, India}
\date{}
\maketitle
\begin{abstract} This paper is about the study of  $F$-transforms based on overlap and grouping maps, residual and co-residual implicator over complete lattice from both constructive and axiomatic approaches. Further,  the duality, basic properties, and the inverse of proposed $F$-transforms have been studied, and axiomatic characterizations of proposed direct $F$-transforms are investigated.
\end{abstract}
\section{Introduction} 
The theory of fuzzy transform ($F$-transform) was firstly introduced by Perfilieva \cite{per}, a notion that piqued the curiosity of many researchers. It has now been greatly expanded upon, and a new chapter in the notion of semi-linear spaces has been opened. The fundamental idea of the $F$-transform is to factorize (or fuzzify) the precise values of independent variables by using a proximity relationship, and to average the precise values of dependent variables to an approximation value (cf., \cite{{per},{irin1}}), from fuzzy sets to parametrized fuzzy sets \cite{st} and from the single variable to the two (or more variables) (cf., \cite{{ma},{mar}, {mar1}, {step}}). Recently, several studies have begun to look into $F$-transforms based on an arbitrary $L$-fuzzy partition of an arbitrary universe (cf., \cite{kh1,mockor1,jir,mock, jiri1,mo,mocko, anan,spt1,trip}), where $L$ is a complete residuated lattice. Among these researches, the concept of a general transformation operator determined by a monadic relation was introduced in \cite{mockor1}, the links between $F$-transforms and {semimodule homomorphisms} were examined in \cite{jir}, while the connections between $F$-transforms and similarity {relations} were discussed in \cite{mocko}. Further, a fascinating relationship of $L$-fuzzy topologies/co-topologies and $L$-fuzzy approximation operators (all of which are ideas employed in the study of an operator-oriented perspective of rough set theory) with $F$-transforms was also discovered in \cite{anan}, while the connection of $L^M$-valued $F$-transforms with $L^M$-valued fuzzy approximation operators and $ML$-graded topologies/co-topologies was discussed in \cite{trip}. Also, the concept of $F$-transforms and $L$-fuzzy pretopologies were examined in \cite{spt1}. In which it has been shown that weaker closure and interior operators, called after \v{C}ech, may also be expressed by using $F$-transforms, implying that $L$-valued $F$-transforms could be utilized in parallel with closure and interior operators as their canonical representation. Also, classes of $F$-transforms taking into account three well-known classes of implicators, namely {$R-,S-,QL-$ implicators} were discussed in \cite{tri}. Several studies in the subject of $F$-transforms applications have been conducted, e.g., trend-cycle estimation \cite{holc}, {data compression \cite{hut}}, numerical solution of partial differential equations \cite{kh}, scheduling \cite{li}, time series \cite{vil}, data analysis \cite{no}, denoising \cite{Ir}, face recognition \cite{roh}, neural network { approaches \cite{ste} and trading} \cite{to}.
\subsection{Motivation of our research}
In contrast to the usual fuzzy logical connectives $t$-norm and $t$-conorms, the overlap and grouping maps can also be regarded as a new structure of classical logic's intersection and union operations on the unit interval. Even though these maps are closely linked to $t$-norm and $t$-conorm, they do not have any nontrivial zero divisors. Recently, several researchers have examined the construction technique and properties of overlap and grouping maps over complete lattices and conducted extensive research. Qiao presented the concepts of overlap and grouping maps over complete lattices in \cite{qiao1} and provided two construction techniques. In \cite{wang}, complete homomorphisms and complete $0_L,1_L$-endomorphisms were used to examine the construction techniques of overlap and grouping maps over complete lattices. Further, the ordinal sums of overlap and grouping maps were discussed in \cite{wang1}. Also, the overlap and grouping maps have been used in various aspects of practical application problems such as in image processing \cite{jurio}, classification \cite{elk}, and decision-making \cite{bus} problems. Specifically, these maps have more advantages than $t$-norm and $t$-conorm in dealing with some real issues. It seems that using the ideas of the overlap and grouping maps in $F$-transform may further open some new areas of application. Accordingly, the study of the theory of $F$-transform using the ideas of such maps is a theme of this paper. 
\subsection{Main contributions}
In this work, we present the theory of $F$-transforms based on overlap and grouping maps, residual and co-residual implicators over complete lattices. Interestingly, under certain conditions, the $F$-transforms introduced in \cite{per,anan,tri} are special cases of proposed $F$-transforms. Further, we study $F$-transforms from constructive and axiomatic approaches based on the above logic operations over complete lattices. The main findings are summarized below: 
\begin{itemize}
\item[$\bullet$] we discuss the duality of the proposed direct $F$-transforms and investigate their basic properties;
\item[$\bullet$] we introduce the inverse of the proposed $F$-transforms and discuss some basic properties; and
\item[$\bullet$] we show a close connection between proposed $F$-transforms and $L$-fuzzy transformation systems and discuss the duality of $L$-fuzzy transformation systems.
\end{itemize}
The remainder of this paper is arranged in the following manner. In Section 2, we recall some key concepts that will be used throughout the main sections. We introduce and examine various classes of direct $F$-transforms determined by overlap and grouping maps over the complete lattice in Section 3. In Section 4, we introduce the inverse of the proposed direct $F$-transforms. In the next section, we characterize proposed direct $F$-transforms from the axiomatic approach. 
\section{Preliminaries}
Herein, we recall the basic ideas related to complete lattices, overlap and grouping maps, $L$-fuzzy sets from \cite{gog,kli,qiao,qiao1,wang,wang1}. Throughout this paper, a complete lattice with the smallest element $0$ and the largest element $1$ is denoted by $L\equiv (L,\vee,\wedge,0,1)$. We start with the following.
\begin{def1}
Let $X$ be a nonempty set. Then an {\bf $L$-fuzzy set} in $X$ is a map $f:X\rightarrow L$.
\end{def1}
The family of all $L$-fuzzy sets in $X$ is denoted by $L^X$. For all $u \in L, \textbf{u}\in L^X$, $\textbf{u}(x) = u, x \in X$ denotes {\bf constant $L$-fuzzy set}. Also, the {\bf core} of an $L$-fuzzy set $f$ is given as a crisp set
$core(f)=\lbrace x\in X,f(x)=1\rbrace.$
If $core(f)\neq \emptyset$, then $f$ is called a {\bf normal $L$-fuzzy set}. For $A\subseteq X$, the {\bf characteristic map} of $A$ is a map $1_A:X\rightarrow \{0,1\}$ such that
\begin{eqnarray*} 
1_A(x)=\begin{cases}
1 &\text{ if } x\in A,\\
0&\text{ otherwise}.
\end{cases} 
\end{eqnarray*}
In the following, we recall and introduce the some basic concepts.
\begin{def1}
An {\bf overlap map} on $L$ is a map ${\theta}:L\times L\rightarrow L$ such that for all $u,v\in L,\{u_i:i\in J\},\{v_i:i\in J\}\subseteq L$
\begin{itemize}
\item[(i)] ${\theta}(u,v)={\theta}(v,u)$,
\item[(ii)] ${\theta}(u,v)=0$ iff $u=0$ or $v=0$,
\item[(iii)] ${\theta}(u,v)=1$ iff $u=1$ and $v=1$,
\item[(iv)] ${\theta}(u,v)\leq {\theta}(u,w)$ if $v\leq w$, and 
\item[(v)] ${\theta}(u,\bigvee\limits_{i\in J}v_i)=\bigvee\limits_{i\in J}{\theta}(u,v_i),{\theta}(\bigwedge\limits_{i\in J}u_i,v)=\bigwedge\limits_{i\in J}{\theta}(u_i,v)$.
\end{itemize}
\end{def1}
If ${\theta}(1,u)=u,\,\forall\,u\in L$, we say that $1$ is a neutral element of ${{\theta}}$. Also, an overlap map is called
\begin{itemize}
\item[(i)] {\bf deflation} if ${\theta}(1,u)\leq u,\,\forall u\in L$,
\item[(ii)] {\bf inflation} if $u\leq{\theta}(1,u),\,\forall u\in L$, and
\item[(iii)] {\bf $EP$-overlap map} if ${\theta}(u,{\theta}(v,w))={\theta}(v,{\theta}(u,w)),\,\forall\,u,v,w\in L$.
\end{itemize}
\begin{exa} (i) Every continuous $t$-norm $\mathcal{T}$ with no nontrivial zero divisors is an overlap map,\\
(ii) ${\theta}_M(u,v)=u\wedge v,\,\forall\,u,v\in L$ on a frame with the prime element $0$ is an overlap map.
\end{exa}
\begin{def1}
A {\bf grouping map} on $L$ is a map ${\eta}:L\times L\rightarrow L$ such that for all $u,v\in L,\{u_i:i\in J\},\{v_i:i\in J\}\subseteq L$
\begin{itemize}
\item[(i)] ${\eta}(u,v)={\eta}(v,u)$,
\item[(ii)] ${\eta}(u,v)=0$ iff $u=0$ and $v=0$,
\item[(iii)] ${\eta}(u,v)=1$ iff $u=1$ or $v=1$,
\item[(iv)] ${\eta}(u,v)\leq {\eta}(u,w)$ if $v\leq w$, and 
\item[(v)] ${\eta}(u,\bigvee\limits_{i\in J}v_i)=\bigvee\limits_{i\in J}{\eta}(u,v_i),{\eta}(\bigwedge\limits_{i\in J}u_i,v)=\bigwedge\limits_{i\in J}{\eta}(u_i,v)$.
\end{itemize}
\end{def1}
If ${\eta}(0,u)=u,\,\forall\,u\in L$, we say that $0$ is a neutral element of $\eta$. Also, a grouping map is called
\begin{itemize}
\item[(i)] {\bf deflation} if ${\eta}(0,u)\geq u,\,\forall u\in L$,
\item[(ii)] {\bf inflation} if $u\geq{\eta}(0,u),\,\forall u\in L$, and
\item[(iii)] {\bf $EP$-grouping map} if ${\eta}(u,{\eta}(v,w))={\eta}(v,{\eta}(u,w)),\,\forall\,u,v,w\in L$.
\end{itemize}
\begin{exa} (i) Every continuous $t$-conorm $\mathcal{S}$ with no nontrivial zero divisors is a grouping map,\\
(ii) ${\eta}_M(u,v)=u\vee v,\,\forall\,u,v\in L$ on a frame with the prime element $0$ is an grouping map.
\end{exa}
\begin{def1}
A {\bf negator} on $L$ is a decreasing map $\mathbf{N}:L\rightarrow L$ such that $\mathbf{N}(0)=1$ and $\mathbf{N}(1)=0$. 
\end{def1}
A negator $\mathbf{N}$ is called {\bf involutive} (strong), if $\mathbf{N}(\mathbf{N}(u))=u,\,\forall\,u\in L$. In addition, a negator $\mathbf{N}$ is called {\bf strict}, if $\mathbf{N}$ is stictly decreasing and continuous, i.e., involutive (as every involutive negator is stictly decreasing and continuous). \\\\
The negator $\mathbf{N}_S(u)=1-u$ on $L=[0,1]$ is usually regarded as the standard negator. For a given negator $\mathbf{N}$, an overlap map ${{\theta}}$ and a grouping map ${\eta}$ are dual with respect to $\mathbf{N}$ if ${\eta}(\mathbf{N}(u),\mathbf{N}(v))=\mathbf{N}({\theta}(u,v)),\theta(\mathbf{N}(u),\mathbf{N}(v))=\mathbf{N}({\eta}(u,v)),\,\forall\,u,v\in L$.
\begin{def1}
Let $\mathbf{N}$ be a negator, ${{\theta}}$ be an overlap map and ${\eta}$ be a grouping map. Then
\begin{itemize}
 \item[(i)] the {\bf residual implicator} induced by an overlap map ${{\theta}}$ is a map $\mathcal{I}_{\theta}:L\times L\rightarrow L$ such that $\mathcal{I}_{\theta}(u,v)=\{w\in L:{\theta}(u,w)\leq v \},\,\forall\,u,v\in L$, and
 \item[(ii)] the {\bf co-residual implicator} induced by a grouping map ${\eta}$ is a map $\mathcal{I}_{\eta}:L\times L\rightarrow L$ such that $\mathcal{I}_{\eta}(u,v)=\{w\in L:{\eta}(u,w)\geq v \},\,\forall\,u,v\in L.$
\end{itemize}
\end{def1} 
\begin{exa} Let $L=[0,1],{\theta}={\theta}_M,{\eta}={\eta}_M$. Then for all $u,v\in L$
\begin{itemize}
\item[(i)] the residual implicator $\mathcal{I}_{{\theta}_M}$ is given as
$
\mathcal{I}_{{\theta}_M}(u,v)=\begin{cases}
1 &\text{ if } u\leq v,\\
v&\text{ otherwise},\,and
\end{cases}$
\item[(ii)] the co-residual implicator $\mathcal{I}_{{\eta}_M}$ is given as
$\mathcal{I}_{{\eta}_M}(u,v)=\begin{cases}
0 &\text{ if } u\leq v,\\
u&\text{ otherwise}.
\end{cases}$
\end{itemize} 
\end{exa}
\begin{lem}
Let ${{\theta}}$ and ${\eta}$ be overlap and grouping maps, respectively. Then ${{\theta}}$ and $\mathcal{I}_{{\theta}}$, ${\eta}$ and $\mathcal{I}_{\eta}$ form two adjoint pairs, respectively, i.e., for all $u,v,w\in L,\,{\theta}(u,v)\leq w\Leftrightarrow u\leq \mathcal{I}_{\theta}(v,w),\,{\eta}(u,v)\geq w\Leftrightarrow u\geq \mathcal{I}_{\eta}(v,w)$, {respectively}.
\end{lem}
\begin{lem}
Let ${{\theta}}$ be an overlap map. Then for all $u,v,w\in L$
\begin{itemize}
\item[(i)] $\mathcal{I}_{\theta}(0,0)=\mathcal{I}_{\theta}(1,1)=1,\mathcal{I}_{\theta}(1,0)=0$,
\item[(ii)] $\mathcal{I}_{\theta}(u,w)\geq \mathcal{I}_{\theta}(v,w),\,\mathcal{I}_{\theta}(w,u)\leq \mathcal{I}_{\theta}(w,v)$ if $u\leq v$,
\item[(iii)] $\mathcal{I}_{{\theta}}$ is an $OP$, $NP$-residual implicator, i.e., $u\leq v\Leftrightarrow \mathcal{I}_{\theta}(u,v)=1, \mathcal{I}_{\theta}(1,u)=u$, respectively iff $1$ is a neutral element of ${{\theta}}$,
\item[(iv)] $\mathcal{I}_{{\theta}}$ is an $IP$-residual implicator, i.e., $\mathcal{I}_{\theta}(u,u)=1$ iff ${{\theta}}$ is a deflation overlap map, 
\item[(v)] $\mathcal{I}_{{\theta}}$ is an $EP$-residual implicator, i.e., $ \mathcal{I}_{\theta}(u,\mathcal{I}_{\theta}(v,w))= \mathcal{I}_{\theta}(v,\mathcal{I}_{\theta}(u,w))$ iff ${{\theta}}$ is an $EP$-overlap map.
\end{itemize}
\end{lem}
\begin{lem} Let ${{\theta}}$ be an overlap map. Then for all $u,v,w\in L,\{u_i:i\in J\},\{v_i:i\in J\}\subseteq L$
\begin{itemize}
\item[(i)] ${\theta}(u,\mathcal{I}_{\theta}(u,v))\leq v,\mathcal{I}_{\theta}(u,{\theta}(u,v))\geq v,\mathcal{I}_{\theta}({\theta}(u,v),0)=\mathcal{I}_{\theta}(u,\mathcal{I}(v,0))$,
\item[(ii)] $ \mathcal{I}_{\theta}(u,\bigwedge\limits_{i\in J}v_i)=\bigwedge\limits_{i\in J}\mathcal{I}_{\theta}(u,v_i), \mathcal{I}_{\theta}(\bigvee\limits_{i\in J}u_i,v)=\bigwedge\limits_{i\in J}\mathcal{I}_{\theta}(u_i,v)$,
\item[(iii)] $ \mathcal{I}_{\theta}(u,\bigvee\limits_{i\in J}v_i)\geq \bigvee\limits_{i\in J}\mathcal{I}_{\theta}(u,v_i)$,
\item[(iv)] ${{\theta}}$ is an $EP$-overlap map iff $\mathcal{I}_{\theta}({\theta}(u,v),w)=\mathcal{I}_{\theta}(u,\mathcal{I}_{\theta}(v,w))$.
\end{itemize}
\end{lem}
If ${{\theta}}$ and ${\eta}$ are dual with respect to an involutive negator $\mathbf{N}$,
then $\mathcal{I}_{{\theta}}$ and $\mathcal{I}_{\eta}$ are dual with respect to the involutive negator $\mathbf{N}$, i.e., $\mathcal{I}_{\eta}(\mathbf{N}(u),\mathbf{N}(v))=\mathbf{N}(\mathcal{I}_{\theta}(u,v))$, $\mathcal{I}_{\theta}(\mathbf{N}(u),\mathbf{N}(v))\linebreak=\mathbf{N}(\mathcal{I}_{\eta}(u,v)),\,\forall\,u,v\in L$. Then we have the following dual properties of $\mathcal{I}_{\eta}$ by the properties of $\mathcal{I}_{{\theta}}$ as follows:
\begin{itemize}
\item[(i)] $\mathcal{I}_{\eta}(0,0)=\mathcal{I}_{\eta}(1,1)=0,\mathcal{I}(0,1)=1$,
\item[(ii)] $\mathcal{I}_{\eta}(u,w)\geq \mathcal{I}_{\eta}(v,w),\,\mathcal{I}_{\eta}(w,u)\leq \mathcal{I}_{\eta}(w,v)$ if $u\leq v$,
\item[(iii)] $\mathcal{I}_{\eta}$ is $OP$ and $NP$-co-residual implicator, i.e., $u\geq v\Leftrightarrow \mathcal{I}_{\eta}(u,v)=0$ and $ \mathcal{I}_{\eta}(0,u)=u$, respectively iff $0$ is a neutral element of $\eta$,
\item[(iv)] $\mathcal{I}_{\eta}$ is an $IP$-co-residual implicator, i.e., $\mathcal{I}_{\eta}(u,u)=0$ iff ${\eta}$ is a deflation grouping map, 
\item[(v)] $\mathcal{I}_{\eta}$ is an $EP$-co-residual implicator, i.e., $ \mathcal{I}_{\eta}(u,\mathcal{I}_{\eta}(v,w))= \mathcal{I}_{\eta}(v,\mathcal{I}_{\eta}(u,w))$ iff ${\eta}$ is an $EP$-grouping map,
\item[(vi)] ${\eta}(u,\mathcal{I}_{\eta}(u,v))\geq v,\mathcal{I}_{\eta}(u,{\eta}(u,v))\leq v,\mathcal{I}_{\eta}({\eta}(u,v),1)=\mathcal{I}_{\eta}(u,\mathcal{I}_{\eta}(v,1))$,
\item[(vii)] $ \mathcal{I}_{\eta}(u,\bigvee\limits_{i\in J}v_i)=\bigvee\limits_{i\in J}\mathcal{I}_{\eta}(u,v_i), \mathcal{I}_{\eta}(\bigwedge\limits_{i\in J}u_i,v)=\bigvee\limits_{i\in J}\mathcal{I}_{\eta}(u_i,v)$,
\item[(viii)] $ \mathcal{I}_{\eta}(u,\bigwedge\limits_{i\in J}v_i)\leq \bigwedge\limits_{i\in J}\mathcal{I}_{\eta}(u,v_i)$,
\item[(ix)] ${\eta}$ is an $EP$-grouping map iff $\mathcal{I}_{\eta}({\eta}(u,v),w)=\mathcal{I}_{\eta}(u,\mathcal{I}_{\eta}(v,w))$.
\end{itemize}
For any $\mathcal{I}_{{\theta}}$ and $\mathcal{I}_{\eta}$, $\mathbf{N}_{\mathcal{I}_{{\theta}}}(u)=\mathcal{I}_{\theta}(u,0)$ and $\mathbf{N}_{\mathcal{I}_{\eta}}(u)=\mathcal{I}_{\eta}(u,1),\forall\, u\in L$ are called the negators induced by $\mathcal{I}_{{\theta}}$ and $\mathcal{I}_{\eta}$, respectively.
Next, we introduce the following notations which are going to be used in subsequent sections.\\\\
{Given an overlap map ${{\theta}}$, a grouping map ${\eta}$, a residual implicator $\mathcal{I}_{{\theta}}$, a co-residual implicator $\mathcal{I}_{\eta}$, {a negator $\mathbf{N}$}, and {$L$-fuzzy} sets $f,g\in L^X$, we define {$L$-fuzzy} sets ${\theta}(f,g),{\eta}(f,g),\mathcal{I}_{\theta}(f,g),\mathcal{I}_{\eta}(f,g)$ and ${\mathbf{N}} (f)$} as follows:
\begin{eqnarray*}
{\theta}(f, g)(x)&=&{\theta}(f(x),g(x)), \forall\,x\in X,\\
{\eta}(f, g)(x)&=&{\eta}(f(x),g(x)), \forall\,x\in X,\\
{\mathcal{I}}_{\theta}(f, g)(x)&=&\mathcal{I}_{\theta}(f(x),g(x)), \forall\,x\in X,\\
{\mathcal{I}}_{\eta}(f, g)(x)&=&\mathcal{I}_{\eta}(f(x),g(x)), \forall\,x\in X,\,\text{and}\\
{({\mathbf{N}} (f))(x)}&=&{\mathbf{N}(f(x)), \forall\,x\in X.}
\end{eqnarray*}
\section{Direct $F$-transforms}
Herein, we consider that ${{\theta}}$ and ${\eta}$ are overlap and grouping maps, and these are dual with respect to an involutive negator $\mathbf{N}$. Also, $\mathcal{I}_{{\theta}}$ and $\mathcal{I}_{\eta}$ are residual and co-residual implicators induced by ${{\theta}}$ and ${\eta}$, respectively, introduced as in Section 2. The main content of this section is to present the concepts of the direct $F$-transforms of $L$-fuzzy sets with respect to the above logic operations. Further, we study and investigate their relationships and discuss their basic properties. We start with the definition of $L$-fuzzy partition from \cite{anan}.
\begin{def1} \label{FP}
A collection $\mathcal{P}$ of normal {$L$-fuzzy} sets $\lbrace A_j:j\in J\rbrace$ is called an {\bf {$L$-fuzzy} partition} of a {nonempty set $X$} if the corresponding collection of ordinary sets $\lbrace core(A_j):j\in J\rbrace$ is partition of $X$. The pair $(X,\mathcal{P})$ is called a {\bf space with {$L$-fuzzy} partition}.
\end{def1}
For an $L$-fuzzy partition $\mathcal{P} = \{A_{j} : j \in J\}$, it is possible to associate the onto index map $k:X\rightarrow J$ such that $k(x)=j$ iff $x\in core(A_{j}).$\\\\
The following is towards the direct $F$-transforms computed with ${{\theta}}$, $\eta$, $\mathcal{I}_{{\theta}}$ and $\mathcal{I}_{\eta}$, where $\mathcal{I}_{{\theta}}$ and $\mathcal{I}_{\eta}$ are residual and co-residual implicators induced by overlap and grouping maps ${\theta},{\eta}$, respectively. Now, we begin with the following.
\begin{def1}\label{UFT}
{Let $\mathcal{P}$ be an {$L$-fuzzy} partition of a set $X$ and $f\in L^X$}. Then
\begin{itemize}
\item[(i)] the {{\bf (direct ${{\theta}}$-upper) $F^{\uparrow,{\theta}}$-transform} of $f$ computed with an overlap map ${{\theta}}$ over the {$L$-fuzzy} partition $\mathcal{P}$} is a collection of lattice elements {$\lbrace F^{\uparrow,{\theta}}_{j}[f]:j\in J\rbrace$} and the $j^{th}$ component of {(direct ${{\theta}}$-upper)} $F^{\uparrow,{\theta}}$-transform is given by
$${F^{\uparrow,{\theta}}_{j}[f]}=\bigvee\limits_{x\in X}{\theta}(A_j(x),f(x)),$$ 
\item[(ii)] the {{\bf (direct $\eta$-lower) $F^{\downarrow,{\eta}}$-transform} of $f$ computed with a grouping map ${\eta}$ over the {$L$-fuzzy} partition $\mathcal{P}$} is a collection of lattice elements {$\lbrace F^{\downarrow,{\eta}}_{j}[f]:j\in J\rbrace$} and the $j^{th}$ component of {(direct $\eta$-lower)} $F^{\downarrow,{\eta}}$-transform is given by
$${F^{\downarrow,{\eta}}_{j}[f]}=\bigwedge\limits_{x\in X}{\eta}(\mathbf{N}(A_j(x)),f(x)),$$
\item[(iii)] the {{\bf (direct $\mathcal{I}_{\eta}$-upper) $F^{\uparrow,\mathcal{I}_{\eta}}$-transform} of $f$ computed with a co-residual implicator $\mathcal{I}_{\eta}$ induced by a grouping map ${\eta}$ over the {$L$-fuzzy} partition $\mathcal{P}$} is a collection of lattice elements {$\lbrace F^{\uparrow,\mathcal{I}_{\eta}}_{j}[f]:j\in J\rbrace$} and the $j^{th}$ component of {(direct $\mathcal{I}_{ \eta}$-upper)} $F^{\uparrow,\mathcal{I}_{\eta}}$-transform is given by
$${F^{\uparrow,\mathcal{I}_{\eta}}_{j}[f]}=\bigvee\limits_{x\in X}{\mathcal{I}_{\eta}}(\mathbf{N}(A_j(x)),f(x)),\, {and }$$ 
\item[(iv)] the {{\bf (direct $\mathcal{I}_{{\theta}}$-lower) $F^{\downarrow,\mathcal{I}_{{\theta}}}$-transform} of $f$ computed with a residual implicator ${\mathcal{I}_{{\theta}}}$ induced by an overlap map ${\theta}$ over the {$L$-fuzzy} partition $\mathcal{P}$} is a collection of lattice elements {$\lbrace F^{\downarrow,\mathcal{I}_{{\theta}}}_{j}[f]:j\in J\rbrace$} and the $j^{th}$ component of {(direct $\mathcal{I}_{{\theta}}$-lower)} $F^{\downarrow,\mathcal{I}_{{\theta}}}$-transform is given by
$${F^{\downarrow,\mathcal{I}_{{\theta}}}_{j}[f]}=\bigwedge\limits_{x\in X}{\mathcal{I}_{{\theta}}}(A_j(x),f(x)).$$
\end{itemize}
\end{def1}
The direct upper $F$-transform computed with a $t$-norm and the direct lower $F$-transform computed with an $R$-implicator proposed in \cite{per,anan,tri} are special cases of $F^{\uparrow,{\theta}}$ and $F^{\downarrow,\mathcal{I}_{{\theta}}}$-transforms, respectively. Also, the direct lower $F$-transform computed with an $S$-implicator proposed in \cite{tri} is a special case of $F^{\downarrow,{\eta}}$-transform. In above-introduced direct $F$-transforms, $F^{\uparrow,\mathcal{I}_{\eta}}$-transform is a new definition. 
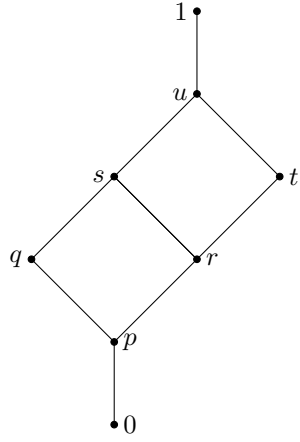
\begin{figure}
 \centering
 \begin{tikzpicture}[scale=.55]
\draw[fill] (1, -4) circle (.08cm) node[right] {$0$};
 \draw[fill] (-1, 0) circle (.08cm) node[left] {$q$};
 \draw[fill] (1, -2) circle (.08cm) node[right] {$p$};
 \draw[fill] (3, 0) circle (.08cm) node[right] {$r$};
 \draw[fill] (1, 2) circle (.08cm) node[left] {$s$};
 \draw[fill] (3, 4) circle (.08cm) node[left] {$u$};
 \draw[fill] (5, 2) circle (.08cm) node[right] {$t$};
 \draw[fill] (3, 6) circle (.08cm) node[left] {$1$};
 \draw (-1,0) -- (1,-2);
 \draw (-1,0) -- (1,2);
 \draw (3,0) -- (1,2);
 \draw (1,-2) -- (3,0);
 \draw (3,0) -- (1,2);
 \draw (1,-2) -- (1,-4);
 \draw (1,2) -- (3,4);
 \draw (3,4) -- (3,6);
 \draw (3,4) -- (5,2);
 \draw (3,0) -- (5,2);
 \end{tikzpicture}
 \caption{Diagram for lattice $L$}
 \label{fig:fig1}
\end{figure}
\begin{exa}\label{exa31} Let $L=\{0,p,q,r,s,t,u,1\}$ be a complete lattice such that $0<p<r<t<u<1,0<p<q<s<u<1$ and $\{q,r\}$ and $\{s,t\}$ are pairwise incomparable (Figure \ref{fig:fig1}). Then $(X,\mathcal{P})$ is a space with an $L$-fuzzy partition $\mathcal{P}$, where $X=\{x_1,x_2,x_3\}$ and $\mathcal{P}=\{A_1,A_2,A_3\}$ such that 
$A_1=\frac{1}{x_1}+\frac{p}{x_2}+\frac{q}{x_3},\,A_2=\frac{s}{x_1}+\frac{1}{x_2}+\frac{u}{x_3},\,A_3=\frac{s}{x_1}+\frac{p}{x_2}+\frac{1}{x_3}.$
Further, let $ f\in L^X$ such that $f=\frac{p}{x_1}+\frac{q}{x_2}+\frac{u}{x_3}$ and $\mathbf{N}$ be an involutive negator such that $\mathbf{N}(0)=1,\mathbf{N}(a)=u,\mathbf{N}(q)=t,\mathbf{N}(r)=s,\mathbf{N}(s)=r,\mathbf{N}(t)=q,\mathbf{N}(u)=p,\mathbf{N}(1)=0$. The the direct $F$-transforms with respect to ${\theta}_M,{\eta}_M,\mathcal{I}_{{\eta}_M},\mathcal{I}_{{\theta}_M}$ are
$F^{\uparrow,{\theta}_M}[f]=\{F^{\uparrow,{\theta}_M}_1[f]=q,
F^{\uparrow,{\theta}_M}_2[f]=u,F^{\uparrow,{\theta}_M}_2[f]=u\},\,
F^{\downarrow,{\eta}_M}[f]=\{F^{\downarrow,{\eta}_M}_1[f]=p,F^{\downarrow,{\eta}_M}_2[f]=r,F^{\downarrow,{\eta}_M}_3[f]=r\},\,
F^{\uparrow,\mathcal{I}_{{\eta}_M}}[f]=\{F^{\uparrow,\mathcal{I}_{{\eta}_M}}_1[f]=u,F^{\uparrow,\mathcal{I}_{{\eta}_M}}_2[f]=r,F^{\uparrow,\mathcal{I}_{{\eta}_M}}_3[f]=u\},\,
F^{\downarrow,\mathcal{I}_{{\theta}_M}}[f]=\{F_1^{\downarrow,\mathcal{I}_{{\theta}_M}}[f]=p,F_2^{\downarrow,\mathcal{I}_{{\theta}_M}}[f]=p,F_3^{\downarrow,\mathcal{I}_{{\theta}_M}}[f]=p\}.$
\end{exa}
\begin{rem}\label{rm32} (i) If $L=[0,1],\,\mathbf{N}=\mathbf{N}_S,{\theta}={\theta}_{M},{\eta}={\eta}_{M},\mathcal{I}_{\eta}=\mathcal{I}_{{\eta}_M}$ and $\mathcal{I}_{\theta}=\mathcal{I}_{{\theta}_M}$, then the $j^{th}$ components of $F^{\uparrow,{\theta}},F^{\downarrow,{\eta}},F^{\uparrow,\mathcal{I}_{\eta}}$ and $F^{\downarrow,\mathcal{I}_{{\theta}}}$-transforms become as follows:
\begin{eqnarray*}
{F^{\uparrow,{\theta}_M}_{j}[f]}&=&\bigvee\limits_{x\in X}(A_j(x)\wedge f(x)),\\
F^{\downarrow,{\eta}_M}_{j}[f]&=&\bigwedge\limits_{x\in X}(( 1-A_j(x))\vee f(x)),\\
F^{\uparrow,\mathcal{I}_{{\eta}_M}}_{j}[f]&=&
\bigvee\limits_{x\in X}\mathcal{I}_{{\eta}_M}((1-A_j(x)), f(x)),\,\text{and}\\
F^{\downarrow,\mathcal{I}_{{\theta}_M}}_{j}[f]&=&\bigwedge\limits_{x\in X}\mathcal{I}_{{\theta}_M}(A_j(x), f(x)),\,\forall\,j\in J,f\in L^X.
\end{eqnarray*}
Obviously $F^{\uparrow,{\theta}_M}$ and $F^{\downarrow,\mathcal{I}_{{\theta}_M}}$-transforms coincide with the special cases of direct upper and lower $F$-transforms proposed in \cite{per,anan,tri}, respectively. Also, $F^{\downarrow,{{\eta}_M}}$-transform coincides with the special case of the direct lower $F$-transform proposed in \cite{tri}.\\\\ 
(ii) If $L=[0,1],{\theta}={\theta}_{M},{\eta}={\eta}_{M},\mathcal{I}_{\eta}=\mathcal{I}_{{\eta}_M}$ and $\mathcal{I}_{\theta}=\mathcal{I}_{{\theta}_M}$, then the $j^{th}$ components of $F^{\uparrow,{\theta}},F^{\downarrow,{\eta}},F^{\uparrow,\mathcal{I}_{\eta}}$ and $F^{\downarrow,\mathcal{I}_{{\theta}}}$-transforms become as follows:
\begin{eqnarray*}
{F^{\uparrow,{\theta}_M}_{j}[f]}&=&\bigvee\limits_{x\in X}(A_j(x)\wedge f(x)),\\
F^{\downarrow,{\eta}_M}_{j}[f]&=&\bigwedge\limits_{x\in X}(\mathbf{N}(A_j(x))\vee f(x)),\\
F^{\uparrow,\mathcal{I}_{{\eta}_M}}_{j}[f]&=&
\bigvee\limits_{x\in X}\mathcal{I}_{{\eta}_M}(\mathbf{N}(A_j(x)), f(x)),\,\text{and}\\
F^{\downarrow,\mathcal{I}_{{\theta}_M}}_{j}[f]&=&\bigwedge\limits_{x\in X}\mathcal{I}_{{\theta}_M}(A_j(x), f(x)),\,\forall\,j\in J,f\in L^X.
\end{eqnarray*}
Obviously $F^{\uparrow,{\theta}_M}$ and $F^{\downarrow,\mathcal{I}_{{\theta}_M}}$-transforms coincide with the special cases of direct upper and lower $F$-transforms proposed in \cite{per,anan,tri}, respectively. Also, $F^{\downarrow,{{\eta}_M}}$-transform coincides with the special case of the direct lower $F$-transform proposed in \cite{tri}.\\\\
(iii) If $L=[0,1],{\theta}=\mathcal{T}$ and ${\eta}=\mathcal{S}$, where $\mathcal{T},\mathcal{S}$ are continuous $t$-norm, $t$-conorm with no nontrivial zero divisors, respectively, then the $j^{th}$ components of $F^{\uparrow,{\theta}},F^{\downarrow,{\eta}},F^{\uparrow,\mathcal{I}_{\eta}}$ and $F^{\downarrow,\mathcal{I}_{{\theta}}}$-transforms become as follows:
\begin{eqnarray*}
F^{\uparrow,\mathcal{T}}_{j}[f]&=&\bigvee\limits_{x\in X}\mathcal{T}(A_j(x),f(x)),\\
{F^{\downarrow,\mathcal{S}}_{j}[f]}&=&\bigwedge\limits_{x\in X}\mathcal{S}(\mathbf{N}( A_j(x)),f(x)),\\
F^{\uparrow,\mathcal{I}_{\mathcal{S}}}_{j}[f]&=&
\bigvee\limits_{x\in X}\mathcal{I}_{\mathcal{S}}(\mathbf{N}(A_j(x)), f(x)),\,\text{and}\\
F^{\downarrow,\mathcal{I}_{\mathcal{T}}}_{j}[f]&=&\bigwedge\limits_{x\in X}\mathcal{I}_{\mathcal{T}}(A_j(x), f(x)),\,\forall\,j\in J,f\in L^X.
\end{eqnarray*}
Obviously $F^{\uparrow,{\mathcal{T}}}$ and $F^{\downarrow,\mathcal{I}_{{\mathcal{T}}}}$-transforms coincide with the direct upper and lower $F$-transforms computed with $t$-norm and $R$-implicator proposed in \cite{per,anan,tri}, respectively. Also, $F^{\downarrow,{{\mathcal{S}}}}$-transform coincide with the direct lower $F$-transform computed with an $S$-implicator proposed in \cite{tri}, respectively.
\end{rem} 
From the above, it is clear that some existing direct $F$-transforms are special cases of the proposed direct $F$-transforms. Among these, some direct $F$-transforms coincide with the proposed direct $F$-transforms and some of the proposed direct $F$-transforms coincide with the special cases of the existing direct $F$-transforms. That is to say; the proposed direct $F$-transforms are more general than other existing ones.
\begin{pro}\label{p32a}
Let ${\theta}$ and ${\eta}$ be dual with respect to an involutive negator $\mathbf{N}$. Then for all $j\in J,f\in L^X$
\begin{itemize}
\item[(i)] $F_j^{\uparrow,{\theta}}[f]=\mathbf{N}(F_j^{\downarrow,{\eta}}[{\mathbf{N}}(f)])$, i.e., $\mathbf{N}(F_j^{\uparrow,{\theta}}[f])=F_j^{\downarrow,{\eta}}[{\mathbf{N}}(f)]$, and
\item[(ii)] $F_j^{\downarrow,{\eta}}[f]=\mathbf{N}(F_j^{\uparrow,{\theta}}[{\mathbf{N}}(f)])$, i.e, $\mathbf{N}(F_j^{\downarrow,{\eta}}[f])=F_j^{\uparrow,{\theta}}[{\mathbf{N}}(f)]$. 
\end{itemize}
\end{pro}
\textbf{Proof:} (i) Let $j\in J$ and $f\in L^X$. Then from Definition {\ref{UFT}}
\begin{eqnarray*}
\mathbf{N}(F_j^{\downarrow,{\eta}}[{\mathbf{N}} (f)])&=&\mathbf{N}(\bigwedge\limits_{x\in X}{\eta}(\mathbf{N}(A_j(x)),({\mathbf{N}} (f))(x)))\\
&=&\mathbf{N}(\bigwedge\limits_{x\in X}{\eta}(\mathbf{N}(A_j(x)),{\mathbf{N}} (f(x))))\\
&=&\bigvee\limits_{x\in X}\mathbf{N}({\eta}(\mathbf{N}(A_j(x)),{\mathbf{N}} (f(x))))\\
&=&\bigvee\limits_{x\in X}{{\theta}}(A_j(x), f(x))\\
&=&F_j^{\uparrow,{\theta}}[f].
\end{eqnarray*} 
Thus $F_j^{\uparrow,{\theta}}[f]=\mathbf{N}(F_j^{\downarrow,{\eta}}[{\mathbf{N}} (f)])$, or that $\mathbf{N}(F_j^{\uparrow,{\theta}}[f])=F_j^{\downarrow,{\eta}}[{\mathbf{N}} (f)]$. Similarly, we can show that $F_j^{\downarrow,{\eta}}[f]=\mathbf{N}(F_j^{\uparrow,{\theta}}[{\mathbf{N}} (f)])$, or that, $\mathbf{N}(F_j^{\downarrow,{\eta}}[f])=F_j^{\uparrow,{\theta}}[{\mathbf{N}} (f)]$.
\begin{pro}\label{p32b}
Let $\mathcal{I}_{{\theta}}$ and $\mathcal{I}_{\eta}$ be dual with respect to an involutive negator $\mathbf{N}$. Then for all $j\in J,f\in L^X$
\begin{itemize}
\item[(i)] $F_j^{\uparrow,\mathcal{I}_{\eta}}[f]=\mathbf{N}(F_j^{\downarrow,\mathcal{I}_{{\theta}}}[{\mathbf{N}}(f)])$, i.e, $\mathbf{N}(F_j^{\uparrow,\mathcal{I}_{\eta}}[f])=F_j^{\downarrow,\mathcal{I}_{{\theta}}}[{\mathbf{N}}(f)]$, and
\item[(ii)] $F_j^{\downarrow,\mathcal{I}_{{\theta}}}[f]=\mathbf{N}(F_j^{\uparrow,\mathcal{I}_{\eta}}[{\mathbf{N}}(f)])$, i.e., $\mathbf{N}(F_j^{\downarrow,\mathcal{I}_{{\theta}}}[f])=F_j^{\uparrow,\mathcal{I}_{\eta}}[{\mathbf{N}}(f)]$. 
\end{itemize}
\end{pro}
\textbf{Proof:} (i) Let $j\in J$ and $f\in L^X$. Then from Definition {\ref{UFT}}
\begin{eqnarray*}
\mathbf{N}(F_j^{\downarrow,\mathcal{I}_{{\theta}}}[{\mathbf{N}} (f)])&=&\mathbf{N}(\bigwedge\limits_{x\in X}\mathcal{I}_{\theta}(A_j(x),({\mathbf{N}} (f))(x)))\\
&=&\mathbf{N}(\bigwedge\limits_{x\in X}\mathcal{I}_{\theta}(A_j(x),{\mathbf{N}} (f(x))))\\
&=&\bigvee\limits_{x\in X}\mathbf{N}(\mathcal{I}_{\theta}(A_j(x),{\mathbf{N}} (f(x))))\\
&=&\bigvee\limits_{x\in X}\mathbf{N}(\mathcal{I}_{\theta}(\mathbf{N}(\mathbf{N}(A_j(x))),{\mathbf{N}} (f(x))))\\
&=&\bigvee\limits_{x\in X}{{\mathcal{I}_{\eta}}}(\mathbf{N}(A_j(x)), f(x))\\
&=&F_j^{\uparrow,{\mathcal{I}_{\eta}}}[f].
\end{eqnarray*} 
Thus $F_j^{\uparrow,{\mathcal{I}_{{\eta}}}}[f]=\mathbf{N}(F_j^{\downarrow,{\mathcal{I}_{\theta}}}[{\mathbf{N}} (f)])$, or that $\mathbf{N}(F_j^{\uparrow,{\mathcal{I}_{{\eta}}}}[f])=F_j^{\downarrow,{\mathcal{I}_{\theta}}}[{\mathbf{N}} (f)]$. Similarly, we can prove that $F_j^{\downarrow,\mathcal{I}_{{\theta}}}[f]=\mathbf{N}(F_j^{\uparrow,\mathcal{I}_{\eta}}[{\mathbf{N}} (f)])$, or that, $\mathbf{N}(F_j^{\downarrow,\mathcal{I}_{{\theta}}}[f])=F_j^{\uparrow,\mathcal{I}_{\eta}}[{\mathbf{N}} (f)]$.\\\\
The above two propositions show that the duality of $F^{\uparrow,{\theta}}_j$ and $F^{\downarrow,{\eta}}_j$, $F^{\uparrow,\mathcal{I}_{\eta}}_j$ and $F^{\uparrow,\mathcal{I}_{{\theta}}}_j$ are
dual with respect $\mathbf{N}$. Generally, duality concept for $F^{\uparrow,{\theta}}_j$
and $F^{\uparrow,\mathcal{I}_{{\theta}}}_j$, $F^{\downarrow,{\eta}}_j$ and $F^{\uparrow,\mathcal{I}_{\eta}}_j$ are not true with respect to $\mathbf{N}$. But they satisfy the following result. For which, we assume $\bigwedge\limits_{b\in L}\mathcal{I}_{\theta}(\mathcal{I}_{\theta}(u,v),v)=u,\,\forall\,u\in L$.
\begin{pro}\label{p32}
Let $\mathbf{N}$ be an involutive neator, $\theta$ and $\eta$ be $EP$-overlap and $EP$-grouping maps, respectively. Then for $j\in J,u\in L,\textbf{u},f\in L^X$
\begin{itemize}
\item[(i)] $ F^{\downarrow,\mathcal{I}_{{\theta}}}_j[f]=\bigwedge\limits_{u\in L}\mathcal{I}_{\theta}(F^{\uparrow,{\theta}}_j[\mathcal{I}_{\theta}(f,\textbf{u})],u),\,F^{\uparrow,{\theta}}_j[f]=\bigwedge\limits_{u\in L}\mathcal{I}_{\theta}(F^{\downarrow,\mathcal{I}_{{\theta}}}_j[\mathcal{I}_{\theta}(f,\textbf{u})],u)$, and
\item[(ii)] $ F^{\downarrow,{\eta}}_j[f]=\bigvee\limits_{u\in L}\mathcal{I}_{\eta}(F^{\uparrow,\mathcal{I}_{\eta}}_j[\mathcal{I}_{\eta}(f,\textbf{u})],u),\,F^{\uparrow,\mathcal{I}_{\eta}}_j[f]=\bigvee\limits_{u\in L}\mathcal{I}_{\eta}(F^{\downarrow,{\eta}}_j[\mathcal{I}_{\eta}(f,\textbf{u})],u)$.
\end{itemize}
\end{pro}
\textbf{Proof:} Let $u\in L$ and $f\in L^X$. Then from Definition \ref{UFT} 
\begin{eqnarray*}
\bigwedge\limits_{u\in L}\mathcal{I}_{\theta}(F^{\uparrow,{\theta}}_j[\mathcal{I}_{\theta}(f,\textbf{u})],u)&=&\bigwedge\limits_{u\in L}\mathcal{I}_{\theta}(\bigvee\limits_{x\in X}{\theta}(A_j(x),\mathcal{I}_{\theta}(f,\textbf{u})(x)),u)\\
&=&\bigwedge\limits_{u\in L}\bigwedge\limits_{x\in X}\mathcal{I}_{\theta}({\theta}(A_j(x),\mathcal{I}_{\theta}(f(x),u)),u)\\
&=&\bigwedge\limits_{u\in L}\bigwedge\limits_{x\in X}\mathcal{I}_{\theta}(A_j(x),\mathcal{I}_{\theta}(\mathcal{I}_{\theta}(f(x),u),u))\\
&=&\bigwedge\limits_{x\in X}\mathcal{I}_{\theta}(A_j(x),\bigwedge\limits_{u\in L}\mathcal{I}_{\theta}(\mathcal{I}_{\theta}(f(x),u),u))\\
&=&\bigwedge\limits_{x\in X}\mathcal{I}_{\theta}(A_j(x),f(x))\\
&=&F^{\downarrow,\mathcal{I}_{{\theta}}}[f].
\end{eqnarray*}
Thus $F^{\downarrow,\mathcal{I}_{{\theta}}}[f]=\bigwedge\limits_{u\in L}\mathcal{I}_{\theta}(F^{\uparrow,{\theta}}_j[\mathcal{I}_{\theta}(f,\textbf{u})],u)$ and 
\begin{eqnarray*}
\bigwedge\limits_{u\in L}\mathcal{I}_{\theta}(F^{\downarrow,\mathcal{I}_{{\theta}}}_j[\mathcal{I}_{\theta}(f,\textbf{u})],u)&=&\bigwedge\limits_{u\in L}\mathcal{I}_{\theta}(\bigwedge\limits_{x\in X}\mathcal{I}_{\theta}(A_j(x),\mathcal{I}_{\theta}(f,\textbf{u})(x)),u)\\
&=&\bigwedge\limits_{u\in L}\mathcal{I}_{\theta}(\bigwedge\limits_{x\in X}\mathcal{I}_{\theta}(A_j(x),\mathcal{I}_{\theta}(f(x),u)),u)\\
&=&\bigwedge\limits_{u\in L}\mathcal{I}_{\theta}(\bigwedge\limits_{x\in X}\mathcal{I}_{\theta}({\theta}(A_j(x),f(x)),u),u)\\
&=&\bigwedge\limits_{u\in L}\mathcal{I}_{\theta}(\mathcal{I}_{\theta}(\bigvee\limits_{x\in X}{\theta}(A_j(x),f(x)),u),u)\\
&=&\bigvee\limits_{x\in X}{\theta}(A_j(x),f(x))\\
&=&F^{\uparrow,{\theta}}[f].
\end{eqnarray*}
Thus $F^{\uparrow,{\theta}}[f]=\bigwedge\limits_{u\in L}\mathcal{I}_{\theta}(F^{\downarrow,\mathcal{I}_{{\theta}}}_j[\mathcal{I}_{\theta}(f,\textbf{u})],u)$.\\\\
(ii) Let $u\in L$ and $f\in L^X$. Then from Definition \ref{UFT} and Propositions \ref{p32a} and \ref{p32b}
\begin{eqnarray*}
F^{\downarrow,{\eta}}[f]&=&\mathbf{N}(F^{\uparrow,{\theta}}[\mathbf{N}(f)])\\
&=&\mathbf{N}(\bigwedge\limits_{u\in L}\mathcal{I}_{\theta}(F^{\downarrow,\mathcal{I}_{{\theta}}}_j[\mathcal{I}_{\theta}(\mathbf{N}(f),\textbf{u})],u))\\
&=&\mathbf{N}(\bigwedge\limits_{u\in L}\mathcal{I}_{\theta}(\bigwedge\limits_{x\in X}\mathcal{I}_{\theta}(A_j(x),\mathcal{I}_{\theta}(\mathbf{N}(f),\textbf{u})(x)),u))\\
&=&\bigvee\limits_{u\in L}\mathbf{N}(\mathcal{I}_{\theta}(\bigwedge\limits_{x\in X}\mathcal{I}_{\theta}(A_j(x),\mathcal{I}_{\theta}(\mathbf{N}(f),\textbf{u})(x)),u))\\
&=&\bigvee\limits_{u\in L}\mathcal{I}_{\eta}(\mathbf{N}(\bigwedge\limits_{x\in X}\mathcal{I}_{\theta}(A_j(x),\mathcal{I}_{\theta}(\mathbf{N}(f),\textbf{u})(x))),\mathbf{N}(u))\\
&=&\bigvee\limits_{u\in L}\mathcal{I}_{\eta}(\bigvee\limits_{x\in X}\mathbf{N}(\mathcal{I}_{\theta}(A_j(x),\mathcal{I}_{\theta}(\mathbf{N}(f),\textbf{u})(x))),\mathbf{N}(u))\\
&=&\bigvee\limits_{u\in L}\mathcal{I}_{\eta}(\bigvee\limits_{x\in X}\mathcal{I}_{\eta}(\mathbf{N}(A_j(x)),\mathbf{N}(\mathcal{I}_{\theta}(\mathbf{N}(f),\textbf{u})(x))),\mathbf{N}(u))\\
&=&\bigvee\limits_{u\in L}\mathcal{I}_{\eta}(\bigvee\limits_{x\in X}\mathcal{I}_{\eta}(\mathbf{N}(A_j(x)),\mathcal{I}_{\eta}(f,\mathbf{N}(\textbf{u}))(x)),\mathbf{N}(u))\\
&=&\bigvee\limits_{u\in L}\mathcal{I}_{\eta}(F^{\uparrow,\mathcal{I}_{\eta}}[\mathcal{I}_{\eta}(f,\mathbf{N}(\textbf{u}))],\mathbf{N}(u)).
\end{eqnarray*}
Thus $ F^{\downarrow,{\eta}}[f]=\bigvee\limits_{u\in L}\mathcal{I}_{\eta}(F^{\uparrow,\mathcal{I}_{\eta}}[\mathcal{I}_{\eta}(f,\mathbf{N}(\textbf{u}))],\mathbf{N}(u))$ and
\begin{eqnarray*}
F^{\downarrow,\mathcal{I}_{\eta}}[f]&=&\mathbf{N}(F^{\downarrow,\mathcal{I}_{{\theta}}}[\mathbf{N}(f)])\\
&=&\mathbf{N}(\bigwedge\limits_{u\in L}\mathcal{I}_{\theta}(F^{\uparrow,{\theta}}_j[\mathcal{I}_{\theta}(\mathbf{N}(f),\textbf{u})],u))\\
&=&\mathbf{N}(\bigwedge\limits_{u\in L}\mathcal{I}_{\theta}(\bigvee\limits_{x\in X}{\theta}(A_j(x),\mathcal{I}_{\theta}(\mathbf{N}(f),\textbf{u})(x)),u))\\
&=&\bigvee\limits_{u\in L}\mathbf{N}(\mathcal{I}_{\theta}(\bigvee\limits_{x\in X}{\theta}(A_j(x),\mathcal{I}_{\theta}(\mathbf{N}(f),\textbf{u})(x)),u))\\
&=&\bigvee\limits_{u\in L}\mathcal{I}_{\eta}(\mathbf{N}(\bigvee\limits_{x\in X}{\theta}(A_j(x),\mathcal{I}_{\theta}(\mathbf{N}(f),\textbf{u})(x))),\mathbf{N}(u))\\
&=&\bigvee\limits_{u\in L}\mathcal{I}_{\eta}(\bigwedge\limits_{x\in X}\mathbf{N}({\theta}(A_j(x),\mathcal{I}_{\theta}(\mathbf{N}(f),\textbf{u})(x))),\mathbf{N}(u))\\
&=&\bigvee\limits_{u\in L}\mathcal{I}_{\eta}(\bigwedge\limits_{x\in X}{\eta}(\mathbf{N}(A_j(x)),\mathbf{N}(\mathcal{I}_{\theta}(\mathbf{N}(f),\textbf{u})(x))),\mathbf{N}(u))\\
&=&\bigvee\limits_{u\in L}\mathcal{I}_{\eta}(\bigwedge\limits_{x\in X}{\eta}(\mathbf{N}(A_j(x)),\mathcal{I}_{\eta}(f,\mathbf{N}(\textbf{u}))(x)),\mathbf{N}(u))\\
&=&\bigvee\limits_{u\in L}{\eta}(F^{\downarrow,{\eta}}[\mathcal{I}_{\eta}(f,\mathbf{N}(\textbf{u}))],\mathbf{N}(u)).
 \end{eqnarray*}
 Thus $F^{\downarrow,\mathcal{I}_{\eta}}[f]=\bigvee\limits_{u\in L}\mathcal{I}_{\eta}(F^{\downarrow,{\eta}}[\mathcal{I}_{\eta}(f,\mathbf{N}(\textbf{u}))],\mathbf{N}(u)).$\\\\
From above three results, we have the following result which present the connections between $F^{\uparrow,{\theta}}$ and $F^{\uparrow,\mathcal{I}_{\eta}}$, $F^{\downarrow,{\eta}}$ and $F^{\downarrow,\mathcal{I}_{{\theta}}}$.
\begin{pro}
Let $\mathbf{N}$ be an involutive negator. Then for $j\in J,u\in L,\textbf{u},f\in L^X$
\begin{itemize}
\item[(i)] $F^{\uparrow,{\theta}}_j[f]=\bigwedge\limits_{u\in L}\mathcal{I}_{\theta}(\mathbf{N}(F^{\uparrow,\mathcal{I}_{\eta}}_j[\mathcal{I}_{\eta}(\mathbf{N}(f),\mathbf{N}(\textbf{u}))]),u)$,
\item[(ii)] $F^{\uparrow,\mathcal{I}_{\eta}}_j[f]=\bigvee\limits_{u\in L}\mathcal{I}_{\eta}(\mathbf{N}(F^{\uparrow,{\theta}}_j[\mathcal{I}_{\eta}(\mathbf{N}(f),\mathbf{N}(\textbf{u}))]),u)$,
\item[(iii)] $ F^{\downarrow,{\eta}}_j[f]=\bigvee\limits_{u\in L}\mathcal{I}_{\eta}(\mathbf{N}(F^{\downarrow,\mathcal{I}_{{\theta}}}_j[\mathcal{I}_{\eta}(\mathbf{N}(f),\mathbf{N}(\textbf{u}))]),u)$, and
\item[(iv)] $ F^{\downarrow,\mathcal{I}_{{\theta}}}_j[f]=\bigwedge\limits_{u\in L}\mathcal{I}_{\theta}(\mathbf{N}(F^{\downarrow,{\eta}}_j[\mathcal{I}_{\eta}(\mathbf{N}(f),\mathbf{N}(\textbf{u}))]),u)$.
\end{itemize}
\end{pro} 
\textbf{Proof:} Propositions \ref{p32a}, \ref{p32b} and \ref{p32} lead to this proof.\\\\
The following are towards the duality of $F^{\uparrow,{\theta}}_j$
and $F^{\uparrow,\mathcal{I}_{{\theta}}}_j$, $F^{\downarrow,{\eta}}_j$ and $F^{\uparrow,\mathcal{I}_{\eta}}_j$ with respect to involutive negators $\mathbf{N}_{\mathcal{I}_{{\theta}}},\mathbf{N}_{\mathcal{I}_{\eta}}$, respectively. 
\begin{pro}\label{p32c}
Let $\mathbf{N}_{\mathcal{I}_{{\theta}}}$ be an involutive negator such that $\mathbf{N}_{\mathcal{I}_{{\theta}}}(.)=\mathcal{I}_{\theta}(.,0)$. Then for all $j\in J,f\in L^X$
\item[(i)] $F_j^{\uparrow,{\theta}}[f]=\mathbf{N}_{\mathcal{I}_{{\theta}}}(F_j^{\downarrow,\mathcal{I}_{{\theta}}}[{\mathbf{N}_{\mathcal{I}_{{\theta}}}}(f)])$, i.e., $\mathbf{N}_{\mathcal{I}_{{\theta}}}(F_j^{\uparrow,{\theta}}[f])=F_j^{\downarrow,\mathcal{I}_{{\theta}}}[{\mathbf{N}_{\mathcal{I}_{{\theta}}}}(f)]$, and
\item[(ii)] $F_j^{\downarrow,\mathcal{I}_{{\theta}}}[f]=\mathbf{N}_{\mathcal{I}_{{\theta}}}(F_j^{\uparrow,{\theta}}[{\mathbf{N}_{\mathcal{I}_{{\theta}}}}(f)])$, i.e, $\mathbf{N}_{\mathcal{I}_{{\theta}}}(F_j^{\downarrow,\mathcal{I}_{{\theta}}}[f])=F_j^{\uparrow,{\theta}}[{\mathbf{N}_{\mathcal{I}_{{\theta}}}}(f)]$. 
\end{pro}
\textbf{Proof:} (i) Let $j\in J$ and $f\in L^X$. Then from Definition {\ref{UFT}}
\begin{eqnarray*}
\mathbf{N}_{\mathcal{I}_{{\theta}}}(F_j^{\downarrow,{\mathcal{I}_{{\theta}}}}[{\mathbf{N}_{\mathcal{I}_{{\theta}}}} (f)])&=&\mathbf{N}_{\mathcal{I}_{{\theta}}}(\bigwedge\limits_{x\in X}{\mathcal{I}_{{\theta}}}(A_j(x),({\mathbf{N}_{\mathcal{I}_{{\theta}}}} (f))(x)))\\
&=&{\mathcal{I}_{{\theta}}}(\bigwedge\limits_{x\in X}{\mathcal{I}_{{\theta}}}(A_j(x),{\mathbf{N}_{\mathcal{I}_{{\theta}}}} (f(x))),0)\\
&=&{\mathcal{I}_{{\theta}}}(\bigwedge\limits_{x\in X}{\mathcal{I}_{{\theta}}}(A_j(x),\mathcal{I}_{\theta}( f(x),0)),0)\\
&=&{\mathcal{I}_{{\theta}}}(\bigwedge\limits_{x\in X}{\mathcal{I}_{{\theta}}}({\theta}(A_j(x),f(x)),0),0)\\
&=&\bigvee\limits_{x\in X}\mathbf{N}_{\mathcal{I}_{{\theta}}}({\mathcal{I}_{{\theta}}}(A_j(x),{\mathbf{N}} (f(x))))\\
&=&\bigvee\limits_{x\in X}{{\theta}}(A_j(x), f(x))\\
&=&F_j^{\uparrow,{\theta}}[f].
\end{eqnarray*} 
Thus $F_j^{\uparrow,{\theta}}[f]=\mathbf{N}_{\mathcal{I}_{{\theta}}}(F_j^{\downarrow,{\mathcal{I}_{{\theta}}}}[{\mathbf{N}_{\mathcal{I}_{{\theta}}}} (f)])$, or that $\mathbf{N}_{\mathcal{I}_{{\theta}}}(F_j^{\uparrow,{\theta}}[f])=F_j^{\downarrow,\mathcal{I}_{{\theta}}}[{\mathbf{N}_{\mathcal{I}_{{\theta}}}} (f)]$.\\\\
(ii) Let $j\in J$ and $f\in L^X$. Then from Definition {\ref{UFT}}
\begin{eqnarray*}
\mathbf{N}_{\mathcal{I}_{{\theta}}}(F_j^{\uparrow,{\mathcal{I}_{{\theta}}}}[{\mathbf{N}_{\mathcal{I}_{{\theta}}}} (f)])&=&\mathbf{N}_{\mathcal{I}_{{\theta}}}(\bigvee\limits_{x\in X}{\theta}(A_j(x),({\mathbf{N}_{\mathcal{I}_{{\theta}}}} (f))(x)))\\
&=&{\mathcal{I}_{{\theta}}}(\bigvee\limits_{x\in X}{\theta}(A_j(x),{\mathbf{N}_{\mathcal{I}_{{\theta}}}} (f(x))),0)\\
&=&\bigwedge\limits_{x\in X}{\mathcal{I}_{{\theta}}}(A_j(x),\mathcal{I}_{\theta}({\mathbf{N}_{\mathcal{I}_{{\theta}}}} (f(x)),0))\\
&=&\bigwedge\limits_{x\in X}{\mathcal{I}_{{\theta}}}(A_j(x),{\mathbf{N}_{\mathcal{I}_{{\theta}}}}({\mathbf{N}_{\mathcal{I}_{{\theta}}}} (f(x))))\\
&=&\bigwedge\limits_{x\in X}{\mathcal{I}_{{\theta}}}(A_j(x), f(x))\\
&=&F_j^{\downarrow,\mathcal{I}_{{\theta}}}[f].
\end{eqnarray*} 
Thus $F_j^{\uparrow,{\theta}}[f]=\mathbf{N}_{\mathcal{I}_{{\theta}}}(F_j^{\downarrow,{\mathcal{I}_{{\theta}}}}[{\mathbf{N}_{\mathcal{I}_{{\theta}}}} (f)])$, or that $\mathbf{N}_{\mathcal{I}_{{\theta}}}(F_j^{\uparrow,{\theta}}[f])=F_j^{\downarrow,\mathcal{I}_{{\theta}}}[{\mathbf{N}_{\mathcal{I}_{{\theta}}}} (f)]$.
\begin{pro}\label{p32d}
Let $\mathbf{N}_{\mathcal{I}_{\eta}}$ be an involutive negator such that $\mathbf{N}_{\mathcal{I}_{\eta}}(.)=\mathcal{I}_{\eta}(.,1)$. Then for all $j\in J,f\in L^X$
\begin{itemize}
\item[(i)] $F_j^{\downarrow,{\eta}}[f]=\mathbf{N}_{\mathcal{I}_{\eta}}(F_j^{\uparrow,\mathcal{I}_{\eta}}[{\mathbf{N}_{\mathcal{I}_{\eta}}}(f)])$, i.e., $\mathbf{N}_{\mathcal{I}_{\eta}}(F_j^{\downarrow,{\eta}}[f])=F_j^{\uparrow,\mathcal{I}_{\eta}}[{\mathbf{N}_{\mathcal{I}_{\eta}}}(f)]$, and
\item[(ii)] $F_j^{\uparrow,\mathcal{I}_{\eta}}[f]=\mathbf{N}_{\mathcal{I}_{\eta}}(F_j^{\downarrow,{\eta}}[{\mathbf{N}_{\mathcal{I}_{\eta}}}(f)])$, i.e, $\mathbf{N}_{\mathcal{I}_{\eta}}(F_j^{\uparrow,\mathcal{I}_{\eta}}[f])=F_j^{\downarrow,{\eta}}[{\mathbf{N}_{\mathcal{I}_{\eta}}}(f)]$.
\end{itemize}
\end{pro}
\textbf{Proof:} Similar to that of Proposition \ref{p32c}.\\\\
Below, we discuss basic results of $F^{\uparrow,{\theta}}_j,F^{\downarrow,{\eta}}_j,F^{\uparrow,\mathcal{I}_{\eta}}_j$ and $F^{\uparrow,\mathcal{I}_{{\theta}}}_j$.
\begin{pro}\label{p31}
Let $\mathcal{P}=\lbrace A_j:j\in J\rbrace,\mathcal{P}'=\lbrace B_{j'}:j'\in J\rbrace$ be $L$-fuzzy partitions of $X$ and $A_j\leq B_{j'},\,\forall \,j,j'\in J$. Then for all $f\in L^X$
\begin{itemize}
\item[(i)] $F^{\uparrow,{\theta}}_j[f]\leq F^{\uparrow,{\theta}}_{j'}[f],F^{\downarrow,{\eta}}_j[f]\geq F^{\downarrow,{\eta}}_{j'}[f]$, and
\item[(ii)] $F^{\uparrow,\mathcal{I}_{\eta}}_j[f]\leq F^{\uparrow,\mathcal{I}_{\eta}}_{j'}[f], F^{\downarrow,\mathcal{I}_{{\theta}}}_j[f]\geq F^{\downarrow,\mathcal{I}_{{\theta}}}_{j'}[f]$. 
\end{itemize}
\end{pro}
\textbf{Proof:} (i) Let $j\in J$ and $f\in L^X$. Then $F^{\uparrow,{\theta}}_j[f]=\bigvee\limits_{x\in X}{\theta}(A_j(x),f(x))\leq \linebreak\bigvee\limits_{x\in X}{\theta}(B_{j'}(x),f(x))=F^{\uparrow,{\theta}}_{j'}[f].$
Thus $F^{\uparrow,{\theta}}_j[f]\leq F^{\uparrow,{\theta}}_{j'}[f]$. Similarly, we can show $F^{\downarrow,{\eta}}_j[f]\geq F^{\downarrow,{\eta}}_{j'}[f]$.\\\\
(ii) Let $j\in J$ and $f\in L^X$. Then 
$F^{\uparrow,\mathcal{I}_{\eta}}_j[f]=\bigvee\limits_{x\in X}\mathcal{I}_{\eta}(\mathbf{N}(A_j(x)),f(x))\leq\linebreak \bigvee\limits_{x\in X}\mathcal{I}_{\eta}(\mathbf{N}(A_j(x)),f(x))=F^{\uparrow,\mathcal{I}_{\eta}}_{j'}[f].$
Thus $F^{\uparrow,\mathcal{I}_{\eta}}_j[f]\leq F^{\downarrow,\mathcal{I}_{\eta}}_{j'}[f]$. Similarly, we can show $F^{\downarrow,\mathcal{I}_{{\theta}}}_j[f]\geq F^{\downarrow,\mathcal{I}_{{\theta}}}_{j'}[f]$.
\begin{pro}\label{p33}
Let $\mathcal{P}$ be an $L$-fuzzy partition of $X$. Then for all $j\in J,f\in L^X,x_j\in core(A_j)$
\begin{itemize}
\item[(i)] $F_j^{\uparrow,{\theta}}[f]\geq {\theta}(1,f(x_j)),F_j^{\downarrow,{\eta}}[f]\leq {\eta}(0,f(x_j))$, and 
\item[(ii)] $F_j^{\uparrow,\mathcal{I}_{\eta}}[f]\geq {\mathcal{I}_{\eta}}(0,f(x_j)),F_j^{\downarrow,\mathcal{I}_{{\theta}}}[f]\leq {\mathcal{I}_{{\theta}}}(1,f(x_j))$.
\end{itemize}
\end{pro} 
\textbf{Proof:} (i) Let $j\in J,x_j\in core(A_j)$ and $f\in L^X$. Then $F_j^{\uparrow,{\theta}}[f]=\bigvee\limits_{x\in X}{\theta}(A_j(x), f(x))\linebreak\geq {\theta}(A_j(x_j), f(x_j))={\theta}(1,f(x_j)).$
Therefore $F_j^{\uparrow,{\theta}}[f]\geq {\theta}(1,f(x_j))$. Similarly, we can show $F_j^{\downarrow,{\eta}}[f]\leq {\eta}(0,f(x_j))$.\\\\
(ii) Let $j\in J,x_j\in core(A_j)$ and $f\in L^X$. Then $F_j^{\uparrow,\mathcal{I}_{\eta}}[f]=\bigvee\limits_{x\in X}\mathcal{I}_{\eta}(\mathbf{N}(A_j(x)), f(x))\linebreak\geq \mathcal{I}_{\eta}(\mathbf{N}(A_j(x_j)), f(x_j))=\mathcal{I}_{\eta}(0,f(x_j)).$
Therefore $F_j^{\uparrow,\mathcal{I}_{\eta}}[f]\geq \mathcal{I}_{\eta}(0,f(x_j))$. Similarly, we can show $F_j^{\downarrow,\mathcal{I}_{{\theta}}}[f]\leq \mathcal{I}_{\theta}(1,f(x_j))$.
\begin{cor}\label{c1} 
Let the conditions of Proposition \ref{p33} be fulfilled and $1,0$ be neutral elements of ${\theta},{\eta}$, respectively. Then for all $f\in L^X,x_j\in core(A_j)$
\begin{itemize}
 \item[(i)] $F_j^{\uparrow,{\theta}}[f]\geq f(x_j),F_j^{\downarrow,{\eta}}[f]\leq f(x_j)$, and 
\item[(ii)] $F_j^{\uparrow,\mathcal{I}_{\eta}}[f]\geq f(x_j),F_j^{\downarrow,\mathcal{I}_{{\theta}}}[f]\leq f(x_j)$.
\end{itemize}
\end{cor}
\textbf{Proof:} Let $1,0$ be neutral elements of ${\theta},{\eta}$, respectively. Then we have ${\theta}(1,f(x_j))=f(x_j),{\eta}(0,f(x_j))=f(x_j),\mathcal{I}_{\eta}(0,f(x_j))=f(x_j)$ and $\mathcal{I}_{\theta}(1,f(x_j))=f(x_j)$. Now, from Proposition \ref{p33}, we have 
\begin{itemize}
 \item[(i)] $F_j^{\uparrow,{\theta}}[f]\geq f(x_j),F_j^{\downarrow,{\eta}}[f]\leq f(x_j)$ {if} $x_j\in core(A_j)$, and
 \item[(ii)] $F_j^{\uparrow,\mathcal{I}_{\eta}}[f]\geq f(x_j), F_j^{\downarrow,\mathcal{I}_{{\theta}}}[f]\leq f(x_j)$ {if} $x_j\in core(A_j)$.
 \end{itemize}
\begin{pro}\label{p34}
Let $\mathcal{P}$ be an $L$-fuzzy partition of $X$. Then for all $j\in J,f,g\in L^X$ and $f\leq g$
\begin{itemize}
\item[(i)] $F_j^{\uparrow,{\theta}}[f]\leq F_j^{\uparrow,{\theta}}[g],F_j^{\downarrow,{\eta}}[f]\leq F_j^{\downarrow,{\eta}}[g]$, and
\item[(ii)] $F_j^{\uparrow,\mathcal{I}_{\eta}}[f]\leq F_j^{\uparrow,\mathcal{I}_{\eta}}[g],F_j^{\downarrow,\mathcal{I}_{{\theta}}}[f]\leq F_j^{\downarrow,\mathcal{I}_{{\theta}}}[g]$.
\end{itemize}
\end{pro}
\textbf{Proof:} (i) Let $j\in J,f,g\in L^X$ and $f\leq g$. Then $F_j^{\uparrow,{\theta}}[f]=\bigvee\limits_{x\in X}{\theta}(A_j(x), f(x))\leq \bigvee\limits_{x\in X}{\theta}(A_j(x), {\eta}(x))=F_j^{\uparrow,{\theta}}[g].$
Thus $F_j^{\uparrow,{\theta}}[f]\leq F_j^{\uparrow,{\theta}}[g]$. Similarly, we can show that $F_j^{\downarrow,{\eta}}[f]\leq F_j^{\downarrow,{\eta}}[g]$.\\\\
(ii) Let $j\in J,f,g\in L^X$ and $f\leq g$. Then $F_j^{\uparrow,\mathcal{I}_{\eta}}[f]=\bigvee\limits_{x\in X}\mathcal{I}_{\eta}(\mathbf{N}(A_j(x)), f(x))\leq \bigvee\limits_{x\in X}\mathcal{I}_{\eta}(\mathbf{N}(A_j(x)), {\eta}(x))=F_j^{\uparrow,\mathcal{I}_{\eta}}[g].$ Thus $F_j^{\uparrow,{\theta}}[f]\leq F_j^{\uparrow,{\theta}}[g]$. Similarly, we can show that $F_j^{\downarrow,\mathcal{I}_{{\theta}}}[f]\leq F_j^{\downarrow,\mathcal{I}_{{\theta}}}[g]$.
\begin{pro}\label{p35}
Let ${{\theta}}$ and ${\eta}$ be $EP$-overlap and $EP$-grouping maps, respectively. Then for all $u\in L,\textbf{u},f\in L^X$
\begin{itemize}
\item[(i)] $F_j^{\uparrow,{\theta}}[{\theta}(\textbf{u},f)]={\theta}({a},F_j^{\uparrow,{\theta}}[f]),F_j^{\downarrow,{\eta}}[{\eta}(\textbf{u},f)]={\eta}({a},F_j^{\downarrow,{\eta}}[f])$, and
\item[(ii)] $F_j^{\uparrow,\mathcal{I}_{\eta}}[\mathcal{I}_{\eta}(\textbf{u},f)]=\mathcal{I}_{\eta}({a},F_j^{\uparrow,\mathcal{I}_{\eta}}[f]),F_j^{\downarrow,\mathcal{I}_{{\theta}}}[\mathcal{I}_{\theta}(\textbf{u},f)]=\mathcal{I}_{\theta}({a},F_j^{\downarrow,\mathcal{I}_{{\theta}}}[f])$.
\end{itemize}
\end{pro}
\textbf{Proof:} (i) Let $\textbf{u},f\in L^X$. Then
\begin{eqnarray*}
F_j^{\uparrow,{\theta}}[{\theta}(\textbf{u},f)]&=&\bigvee\limits_{x\in X}{\theta}(A_j(x), {\theta}(\textbf{u},f)(x))=\bigvee\limits_{x\in X}{\theta}(A_j(x), {\theta}(u,f(x)))\\
&=&\bigvee\limits_{x\in X}{\theta}(u, {\theta}(A_j(x),f(x)))={\theta}(u,\bigvee\limits_{x\in X} {\theta}(A_j(x),f(x)))\\
&=& {\theta}({u},F_j^{\uparrow,{\theta}}[f]).
\end{eqnarray*}
Therefore $F_j^{\uparrow,{\theta}}[{\theta}(\textbf{u},f)]={\theta}({u},F_j^{\uparrow,{\theta}}[f])$. Similarly, we can show $F_j^{\downarrow,{\eta}}[{\eta}(\textbf{u},f)]= {\eta}({u},F_j^{\downarrow,{\eta}}[f])$.\\\\
(ii) Let $u\in L$ and $\textbf{u},f\in L^X$. Then 
\begin{eqnarray*}
F_j^{\uparrow,\mathcal{I}_{\eta}}[f]&=&\bigvee\limits_{x\in X}\mathcal{I}_{\eta}(\mathbf{N}(A_j(x)), \mathcal{I}_{\eta}(\textbf{u},f)(x))=\bigvee\limits_{x\in X}\mathcal{I}_{\eta}(\mathbf{N}(A_j(x)), \mathcal{I}_{\eta}(u,f(x)))\\
&=&\bigvee\limits_{x\in X}\mathcal{I}_{\eta}(u, \mathcal{I}_{\eta}(\mathbf{N}(A_j(x)),f(x)))=\mathcal{I}_{\eta}(u,\bigvee\limits_{x\in X} \mathcal{I}_{\eta}(\mathbf{N}(A_j(x)),f(x)))\\
&=& \mathcal{I}_{\eta}({u},F_j^{\uparrow,\mathcal{I}_{\eta}}[f]).
\end{eqnarray*}
Therefore $F_j^{\uparrow,\mathcal{I}_{\eta}}[\mathcal{I}_{\eta}(\textbf{u},f)]=\mathcal{I}_{\eta}({u},F_j^{\uparrow,\mathcal{I}_{\eta}}[f])$. Similarly, we can show
$F_j^{\downarrow,\mathcal{I}_{{\theta}}}[\mathcal{I}_{\theta}(\textbf{u},f)]= \mathcal{I}_{\theta}({u},F_j^{\downarrow,\mathcal{I}_{{\theta}}}[f])$.
\begin{pro}\label{p36}
Let $\mathcal{P}$ be an $L$-fuzzy partition of $X$. Then for all $\textbf{u},f\in L^X,\{f_j:j\in J\}\subseteq L^X$
\begin{itemize}
\item[(i)] $F_j^{\uparrow,{\theta}}[\bigvee\limits_{k\in J}f_k]=\bigvee\limits_{k\in J}F_j^{\uparrow,{\theta}}[f_k],F_j^{\downarrow,{\eta}}[\bigwedge\limits_{k\in J}f_k]=\bigwedge\limits_{k\in J}F_j^{\downarrow,{\eta}}[f_k]$, and
\item[(ii)] $F_j^{\uparrow,\mathcal{I}_{\eta}}[\bigvee\limits_{k\in J}f_k]=\bigvee\limits_{k\in J}F_j^{\uparrow,\mathcal{I}_{\eta}}[f_k],F_j^{\downarrow,\mathcal{I}_{{\theta}}}[\bigwedge\limits_{k\in J}f_k]=\bigwedge\limits_{k\in J}F_j^{\downarrow,\mathcal{I}_{{\theta}}}[f_k]$.
\end{itemize}
\end{pro}
 \textbf{Proof:} (i) Let $\{f_k:k\in J\}\subseteq L^X$. Then $
F_j^{\uparrow,{\theta}}[\bigvee\limits_{k\in }f_k]=\bigvee\limits_{x\in X}{\theta}(A_j(x), \bigvee\limits_{k\in J}f_k(x))=\bigvee\limits_{x\in X}\bigvee\limits_{k\in J}{\theta}(A_j(x), f_k(x))=\bigvee\limits_{k\in J} F_k^{\uparrow,{\theta}}[f_k].$
Therefore $F_j^{\uparrow,{\theta}}[\bigvee\limits_{k\in }f_k]=\bigvee\limits_{k\in J} F_k^{\uparrow,{\theta}}[f_k]$. Similarly, we obtain $F_j^{\downarrow,{\eta}}[\bigwedge\limits_{j\in J}f_k]= \bigwedge\limits_{j\in J}F_j^{\downarrow,{\eta}}[f_k]$.\\\\
(ii) Let $\{f_k:k\in J\}\subseteq L^X$. Then
$F_j^{\uparrow,\mathcal{I}_{\eta}}[\bigvee\limits_{k\in }f_k]=\bigvee\limits_{x\in X}\mathcal{I}_{\eta}(\mathbf{N}(A_j(x)), \bigvee\limits_{k\in J}f_k(x))=\bigvee\limits_{x\in X}\bigvee\limits_{k\in J}\mathcal{I}_{\eta}(\mathbf{N}(A_j(x)), f_k(x))=\bigvee\limits_{k\in J} F_k^{\uparrow,\mathcal{I}_{\eta}}[f_k].$
Therefore $F_j^{\uparrow,\mathcal{I}_{\eta}}[\bigvee\limits_{k\in }f_k]=\bigvee\limits_{k\in J} F_k^{\uparrow,\mathcal{I}_{\eta}}[f_k]$. Similarly, we obtain $F_j^{\downarrow,\mathcal{I}_{{\theta}}}[\bigwedge\limits_{j\in J}f_k]= \bigwedge\limits_{j\in J}F_j^{\downarrow,\mathcal{I}_{{\theta}}}[f_k]$.
\begin{pro}\label{p313}
Let $\mathcal{P}$ be an $L$-fuzzy partition of $X$. Then for $u\in L,\textbf{u}\in L^X$
\begin{itemize}
\item[(i)] $F_j^{\uparrow,{\theta}}[\textbf{u}]={\theta}(1,u),F_j^{\downarrow,\mathcal{I}_{{\theta}}}[\textbf{u}]=\mathcal{I}_{\theta}(1,u)$, and
\item[(ii)] $F_j^{\downarrow,{\eta}}[\textbf{u}]={\eta}(\bigwedge\limits_{x\in X}(\mathbf{N}(A_j(x)),u),F_j^{\uparrow,\mathcal{I}_{\eta}}[\textbf{u}]=\mathcal{I}_{\eta}(\bigwedge\limits_{x\in X}(\mathbf{N}(A_j(x)),u)$. In addition, for a strict negator $\mathbf{N}$, $F_j^{\downarrow,{\eta}}[\textbf{u}]={\eta}(0,u),F_j^{\downarrow,\mathcal{I}_{\eta}}[\textbf{u}]=\mathcal{I}_{\eta}(0,u)$. \end{itemize}
\end{pro}
\textbf{Proof:} (i) Let $u\in L$ and $\textbf{u}\in L^X$. Then
$F_j^{\uparrow,{\theta}}[\textbf{u}]=\bigvee\limits_{x\in X}{\theta}(A_j(x), \textbf{u}(x))={\theta}(\bigvee\limits_{x\in X}A_j(x), {a})\linebreak={\theta}(1,u).$
Thus $F_j^{\uparrow,{\theta}}[\textbf{u}]={\theta}(1,u)$. Similarly, we can show $F_j^{\downarrow,\mathcal{I}_{{\theta}}}[\textbf{u}]=\mathcal{I}_{\theta}(1,u)$.\\\\
(ii) Let $u\in L$ and $\textbf{u}\in L^X$. Then $
F_j^{\downarrow,{\eta}}[\textbf{u}]=\bigwedge\limits_{x\in X}{\eta}(\mathbf{N}(A_j(x)), \textbf{u}(x))={\eta}(\bigwedge\limits_{x\in X}\mathbf{N}(A_j(x)), {a}).$
Thus $F_j^{\downarrow,{\eta}}[\textbf{u}]={\eta}(\bigwedge\limits_{x\in X}\mathbf{N}(A_j(x)), {a})$. Now, let $\mathbf{N}$ be a strict negator. Then we obtain $F_j^{\downarrow,{\eta}}[\textbf{u}]={\eta}(\bigwedge\limits_{x\in X}\mathbf{N}(A_j(x)), {a})={\eta}(\mathbf{N}(\bigvee\limits_{x\in X}A_j(x)), {a})= {\eta}(\mathbf{N}(1), {a})={\eta}(0,u).$
Thus $F_j^{\downarrow,{\eta}}[\textbf{u}]={\eta}(0,u)$. Similarly, we can show that $F_j^{\uparrow,\mathcal{I}_{\eta}}[\textbf{u}]=\mathcal{I}_{\eta}(\bigwedge\limits_{x\in X}\mathbf{N}(A_j(x)), {a})$ and for a strict negator $\mathbf{N}$, $\mathcal{I}_{\eta}[\textbf{u}]=\mathcal{I}_{\eta}(0,u)$.
\begin{cor} Let the conditions of Proposition \ref{p313} be fulfilled and $1,0$ be neutral elements of ${\theta},{\eta}$, respectively. Then for all $f\in L^X$
\begin{itemize}
\item[(i)] $F_j^{\uparrow,{\theta}}[\textbf{u}]=u,F_j^{\downarrow,{\eta}}[\textbf{u}]=u$, and
\item[(ii)] $F_j^{\uparrow,\mathcal{I}_{\eta}}[\textbf{u}]=u,F_j^{\downarrow,\mathcal{I}_{{\theta}}}[\textbf{u}]=u$.
\end{itemize}
\end{cor}
\textbf{Proof:} Let $1,0$ be neutral elements of ${\theta},{\eta}$, respectively. Then we have, ${\theta}(1,u)=u,{\eta}(0,u)=u,\mathcal{I}_{\eta}(0,u)=u$ and $\mathcal{I}_{\theta}(1,u)=u$. Also, from Proposition \ref{p313}, we have 
\begin{itemize}
 \item[(i)] $F_j^{\uparrow,{\theta}}[\textbf{u}]=u,F_j^{\downarrow,{\eta}}[\textbf{u}]=u$, and
\item[(ii)] $F_j^{\uparrow,\mathcal{I}_{\eta}}[\textbf{u}]=u,F_j^{\downarrow,\mathcal{I}_{{\theta}}}[\textbf{u}]=u$.
 \end{itemize}
 From Proposition \ref{p313}, we have the follwoing.
\begin{pro}\label{p355}
Let $\mathcal{P}$ be an $L$-fuzzy partition of $X$. Then for all $u\in L,\textbf{u}\in L^X$
\begin{itemize}
\item[(i)] $F_j^{\downarrow,{\eta}}[\textbf{u}]={\eta}(0,u)$ iff $F_j^{\downarrow,{\eta}}[0_X]=0$, and
\item[(ii)]$F_j^{\uparrow,\mathcal{I}_{\eta}}[\textbf{u}]=\mathcal{I}_{\eta}(0,u)$ iff $F_j^{\uparrow,\mathcal{I}_{\eta}}[1_X]=1$.
\end{itemize}
\end{pro}
\textbf{Proof:} (i) Let $F_j^{\downarrow,{\eta}}[\textbf{u}]={\eta}(0,u),\,\forall\,u\in L,\textbf{u}\in L^X$. Then by assuming $\textbf{u}=0_X$, we have $F_j^{\downarrow,{\eta}}[0_X]={\eta}(0,0)=0$. Thus $F_j^{\downarrow,{\eta}}[0_X]=0$. Conversely, from Proposition \ref{p34}(ii), we have $F_j^{\downarrow,{\eta}}[0_X]=0\Leftrightarrow {\eta}(\bigwedge\limits_{x\in X}(\mathbf{N}(A_j(x)),0)=0\Leftrightarrow\bigwedge\limits_{x\in X}(\mathbf{N}(A_j(x))=0.$
Therefore $F_j^{\downarrow,{\eta}}[\textbf{u}]= {\eta}(\bigwedge\limits_{x\in X}(\mathbf{N}(A_j(x)),u)={\eta}(0,u)$. Thus $F_j^{\downarrow,{\eta}}[\textbf{u}]={\eta}(0,u)$. \\\\
(ii) Let $F_j^{\uparrow,\mathcal{I}_{\eta}}[\textbf{u}]=\mathcal{I}_{\eta}(0,u),\,\forall\,u\in L,\textbf{u}\in L^X$. Then by assuming $\textbf{u}=1_X$, we have $F_j^{\uparrow,\mathcal{I}_{\eta}}[1_X]=\mathcal{I}_{\eta}(0,1)=1$. Thus $F_j^{\uparrow,\mathcal{I}_{\eta}}[1_X]=1$. Conversely, from Proposition \ref{p34}(ii), we have
$F_j^{\uparrow,\mathcal{I}_{\eta}}[1_X]=1\Leftrightarrow \mathcal{I}_{\eta}(\bigwedge\limits_{x\in X}(\mathbf{N}(A_j(x)),1)=1\Leftrightarrow\bigwedge\limits_{x\in X}(\mathbf{N}(A_j(x))=0.$
Therefore $F_j^{\downarrow,{\eta}}[\textbf{u}]= \mathcal{I}_{\eta}(\bigwedge\limits_{x\in X}(\mathbf{N}(A_j(x)),u)=\mathcal{I}_{\eta}(0,u)$. Thus $F_j^{\downarrow,{\eta}}[\textbf{u}]=\mathcal{I}_{\eta}(0,u)$. \\\\
The following results are towards the characterization of the components of the direct $F$-transforms of an original $L$-fuzzy set as its lower and upper mean values give the greatest and the least elements to certain sets, respectively.
\begin{pro}\label{p366}
 Let {$\mathcal{P}$} be an {$L$-fuzzy} partition of $X$ and $f\in L^X$. Then 
 \begin{itemize}
\item[(i)] the $j^{th}$ component of $F^{\uparrow,{\theta}}$-transform of $f$ is the least element of the set 
$U_j=\lbrace u\in L:{\theta}(A_j(x),f(x))\leq u,\,\forall x\in X\rbrace,\,j\in J$, and 
\item[(ii)] the $j^{th}$ component of $F^{\downarrow,{\eta}}$-transform of $f$ is the greatest element of the set 
$V_j=\lbrace v\in L:v\leq{\eta}(\mathbf{N}(A_j(x)),f(x)),\,\forall x\in X\rbrace,\,j\in J.$
\end{itemize}
\end{pro}
\textbf{Proof:} (i) To prove this, we need to show that $F_j^{\uparrow,{\theta}}[f]\in U_j$ and $F_j^{\uparrow,{\theta}}[f]\leq u$. It follows from Definition {\ref{UFT}(i)} that
$F_j^{\uparrow,{\theta}}[f]=\bigvee\limits_{x\in X}{\theta}(A_j(x),f(x))\geq {\theta}(A_j(x),f(x)).$
Thus $F_j^{\uparrow,{\theta}}[f]\in U_j$. Now, let $u\in L,x\in X$. Then from the given condition
${\theta}(A_j(x),f(x))\leq u\Rightarrow \bigvee\limits_{x\in X}{\theta}(A_j(x),f(x))\leq u\Rightarrow F_j^{\uparrow,{\theta}}[f]\leq u.$
Thus the $j^{th}$ component of $F^{\uparrow,{\theta}}$-transform is the least element of the set $U_j$.\\\\
(ii) To prove this, we need to show that $F_j^{\downarrow,{\eta}}[f]\in V_j$ and $v\leq F_j^{\downarrow,{\eta}}[f]$. It follows from Definition {\ref{UFT}(ii)} that
$F_j^{\downarrow,{\eta}}[f]=\bigwedge\limits_{x\in X}{\eta}(\mathbf{N}(A_j(x)),f(x))\leq {\eta}(\mathbf{N}(A_j(x)),f(x)).$
Thus $F_j^{\downarrow,{\eta}}[f]\in V_j$. Now, let $v\in L,x\in X$. Then from the given condition
$v\leq {\eta}(\mathbf{N}(A_j(x)),f(x))\Rightarrow v\leq\bigwedge\limits_{x\in X}{\eta}(\mathbf{N}(A_j(x)),f(x))\Rightarrow v\leq F_j^{\downarrow,{\eta}}[f].$
Thus the $j^{th}$ component of $F^{\downarrow,{\eta}}$-transform is the greatest element of the set $V_j$.
\begin{pro}\label{p377}
 Let {$\mathcal{P}$} be an {$L$-fuzzy} partition of $X$ and $f\in L^X$. Then 
 \begin{itemize}
\item[(i)] the $j^{th}$ component of $F^{\uparrow,\mathcal{I}_{\eta}}$-transform of $f$ is the least element of the set 
$U_j=\lbrace u\in L:{\mathcal{I}_{\eta}}(\mathbf{N}(A_j(x)),f(x))\leq u,\,\forall x\in X\rbrace,\,j\in J$, and
\item[(ii)] the $j^{th}$ component of $F^{\downarrow,\mathcal{I}_{{\theta}}}$-transform of $f$ is the greatest element of the set 
$V_j=\lbrace v\in L:v\leq{\mathcal{I}_{{\theta}}}(A_j(x),f(x)),\,\forall x\in X\rbrace,\,j\in J.$
\end{itemize}
\end{pro}
\textbf{Proof:} (i) To prove this, we need to show that $F_j^{\uparrow,\mathcal{I}_{\eta}}[f]\in U_j$ and $F_j^{\uparrow,\mathcal{I}_{\eta}}[f]\leq u$. It follows from Definition {\ref{UFT}(i)} that $
F_j^{\uparrow,\mathcal{I}_{\eta}}[f]=\bigvee\limits_{x\in X}{\mathcal{I}_{\eta}}(\mathbf{N}(A_j(x)),f(x))\geq {\mathcal{I}_{\eta}}(\mathbf{N}(A_j(x)),f(x)).$ Thus $F_j^{\uparrow,\mathcal{I}_{\eta}}[f]\in U_j$. Now, let $u\in L,x\in X$. Then from the given condition
${\mathcal{I}_{\eta}}(\mathbf{N}(A_j(x)),f(x))\leq u\Rightarrow \bigvee\limits_{x\in X}{\mathcal{I}_{\eta}}(\mathbf{N}(A_j(x)),f(x))\leq u\Rightarrow F_j^{\uparrow,\mathcal{I}_{\eta}}[f]\leq u.$
Thus the $j^{th}$ component of $F^{\uparrow,\mathcal{I}_{\eta}}$-transform is the least element of the set $U_j$.\\\\
(ii) To prove this, we need to show that $F_j^{\downarrow,\mathcal{I}_{{\theta}}}[f]\in V_j$ and $v\leq F_j^{\downarrow,\mathcal{I}_{{\theta}}}[f]$. It follows from Definition {\ref{UFT}(ii)} that
$F_j^{\downarrow,\mathcal{I}_{{\theta}}}[f]=\bigwedge\limits_{x\in X}{\mathcal{I}_{{\theta}}}(A_j(x),f(x))\leq {\mathcal{I}_{{\theta}}}(A_j(x),f(x)).$
Thus $F_j^{\downarrow,\mathcal{I}_{{\theta}}}[f]\in V_j$. Now, let $v\in L,x\in X$. Then from the given condition
$v\leq {\mathcal{I}_{{\theta}}}(A_j(x),f(x))\Rightarrow v\leq\bigwedge\limits_{x\in X}{\mathcal{I}_{{\theta}}}(A_j(x),f(x))\Rightarrow v\leq F_j^{\downarrow,\mathcal{I}_{{\theta}}}[f].$
Thus the $j^{th}$ component of $F^{\downarrow,\mathcal{I}_{{\theta}}}$-transform is the greatest element of the set $V_j$.
\begin{pro}
 Let conditions of Proposition \ref{p366} be fulfilled, {${{\theta}}$ and ${\eta}$} be deflation overlap and deflation grouping maps, respectively. Then for all $u\in U_j,v\in V_j$
 \begin{itemize}
\item[(i)] $\bigwedge\limits_{x\in X}\mathcal{I}_{\theta}({\theta}(A_j(x),f(x)),u))=1$
and $j^{th}$ component of $F^{\uparrow,{\theta}}$-transform is the smallest such $u$, and
\item[(ii)] $\bigwedge\limits_{x\in X}{\mathcal{I}_{\eta}}({\eta}(\mathbf{N}(A_j(x)),f(x)),v)=0$ and $j^{th}$ component of $F^{\downarrow,{\eta}}$-transform is the greatest such $v$.
\end{itemize}
 \end{pro}
 \textbf{Proof:} (i) Let $j\in J$. Then for all $x\in X$, ${\theta}(A_j(x),f(x))\leq u$, or that, \linebreak$\bigwedge\limits_{x\in X}{\mathcal{I}_{{\theta}}}({\theta}(A_j(x),f(x)),u)=1$, as ${\mathcal{I}_{{\theta}}}$ is an $IP$-residual implicator.\\\\
 (ii) Let $j\in J$. Then for all $x\in X$, $ {\eta}(\mathbf{N}(A_j(x)),f(x))\geq v$, or that,\linebreak $\bigwedge\limits_{x\in X}{\mathcal{I}_{\eta}}({\eta}(\mathbf{N}(A_j(x)),f(x)),v)=0$, as ${\mathcal{I}_{\eta}}$ is an $IP$-co-residual implicator.
\begin{pro}\label{3.4}
{Let conditions of Proposition \ref{p377} be fulfilled, {${{\theta}}$ and ${\eta}$} be deflation overlap and deflation grouping maps, respectively. Then for all $u\in U_j,v\in V_j$ }
 \begin{itemize}
 \item[(i)] $\bigwedge\limits_{x\in X}{\mathcal{I}_{\eta}}(u,{\mathcal{I}_{\eta}}(\mathbf{N}(A_j(x)),f(x)))=0$ and $j^{th}$ component of $F^{\uparrow,\mathcal{I}_{\eta}}$-transform is the smallest such $u$, and
\item[(ii)] $\bigwedge\limits_{x\in X}{\mathcal{I}_{{\theta}}}(v,{\mathcal{I}_{{\theta}}}(A_j(x),f(x)))=1$ and $j^{th}$ component of $F^{\downarrow,\mathcal{I}_{{\theta}}}$-transform is the greatest such $v$.
\end{itemize}
 \end{pro}
 \textbf{Proof:} (i) Let $j\in J$. Then for all $x\in X$, ${\mathcal{I}_{\eta}}(\mathbf{N}(A_j(x)),f(x))\leq u$, or that, \linebreak $\bigwedge\limits_{x\in X}{\mathcal{I}_{\eta}}(u,{\mathcal{I}_{\eta}}(\mathbf{N}(A_j(x)),f(x)))=0$, as ${\mathcal{I}_{\eta}}$ is an $IP$-co-residual implicator.\\\\
 (ii) Let $j\in J$. Then for all $x\in X$, $ {\mathcal{I}_{{\theta}}}(A_j(x),f(x))\geq v$, or that,\linebreak $\bigwedge\limits_{x\in X}{\mathcal{I}_{{\theta}}}(v,{\mathcal{I}_{{\theta}}}(A_j(x),f(x)))=1$, as ${\mathcal{I}_{{\theta}}}$ is an $IP$-residual implicator.
\section{Inverse $F$-transforms}
In this section, we introduce the concepts of the inverse $F$-transforms computed with overlap and grouping maps, residual and co-residual implicators over $L$, respectively. Further, we discuss their properties. Now, we begin with the following.
\begin{def1}\label{IT}{
Let $(X,\mathcal{P})$ be a space with an $L$-fuzzy partition $\mathcal{P}$ and $f\in L^X$, where $\mathcal{P}=\{A_j\in L^X:j\in J\}$. Further, let $F^{\uparrow,{\theta}}_j [f]$ and $F^{\downarrow,{\eta}}_j [f]$ be the $j^{th}$ components of $F^{\uparrow,{\theta}}$-transform of $f$ computed with an overlap map ${{\theta}}$ over $\mathcal{P}$ and $F^{\downarrow,{\eta}}$-transform
of $f$ computed with a grouping map ${\eta}$ over $\mathcal{P}$, respectively. Then
\begin{itemize}
\item[(i)] the {\bf inverse (upper) $F^{\uparrow,{\theta}}$-transform} of $f$ computed with a residual implication $\mathcal{I}_{{\theta}}$ over a fuzzy partition $\mathcal{P}$ is a map {$\hat{f}^{\uparrow,\mathcal{I}_{{\theta}}}:L^X\rightarrow L^X$} such that
$${\hat{f}^{\uparrow,\mathcal{I}_{{\theta}}}(x)}=\bigwedge\limits_{j\in J}{\mathcal{I}_{{\theta}}}(A_j(x),F^{\uparrow,{\theta}}_j[f]),$$
\item[(ii)] the {\bf inverse (lower) $F^{\downarrow,\mathcal{I}_{{\theta}}}$-transform} of $f$ computed with an overlap map ${{\theta}}$ over a fuzzy partition $\mathcal{P}$ is a map {$\hat{f}^{\downarrow,{\theta}}:L^X\rightarrow L^X$} such that
$${\hat{f}^{\downarrow,{\theta}}(x)}=\bigvee\limits_{j\in J}{\theta}(A_j(x),F^{\downarrow,\mathcal{I}_{{\theta}}}_j[f]),$$
\item[(iii)] the {\bf inverse (upper) $F^{\uparrow,\mathcal{I}_{\eta}}$-transform} of $f$ computed with a grouping map ${\eta}$ over a fuzzy partition $\mathcal{P}$ is a map {$\hat{f}^{\uparrow,{\eta}}:L^X\rightarrow L^X$} such that
$${\hat{f}^{\uparrow,{\eta}}}=\bigwedge\limits_{j\in J}{\eta}(\mathbf{N}(A_j(x)),F^{\uparrow,\mathcal{I}_{\eta}}_j[f]),\,{and}$$
\item[(iv)] the {\bf inverse (lower) $F^{\downarrow,{\eta}}$-transform} of $f$ computed with a co-residual implicator $\mathcal{I}_{\eta}$ over
a fuzzy partition $\mathcal{P}$ is a map {$\hat{f}^{\downarrow,\mathcal{I}_{\eta}}:L^X\rightarrow L^X$} such that
$${\hat{f}^{\downarrow,\mathcal{I}_{\eta}}(x)}=\bigvee\limits_{j\in J}{\mathcal{I}_{\eta}}(\mathbf{N}(A_j(x)),F^{\downarrow,{\eta}}_j[f]).$$
\end{itemize}}
\end{def1} 
The inverse $F$-transforms computed with a $t$-norm and an $R$-implicator proposed in \cite{per,anan,tri} are special cases of the proposed inverse $F$-transforms with respect to ${{\theta}}$ and $\mathcal{I}_{{\theta}}$. In the above-introduced inverse $F$-transforms, $\hat{f}^{\downarrow,{\eta}}$, $\hat{f}^{\downarrow,\mathcal{I}_{\eta}}$ are new definitions.
\begin{exa}\label{exa41}
In continuation to Example \ref{exa31}, the inverse $F$-transforms with respect to $\mathcal{I}_{{\theta}_M},{\theta}_M,{\eta}_M,\mathcal{I}_{{\eta}_M}$ are
$\hat{f}^{\uparrow,\mathcal{I}_{{\theta}_M}}=\frac{q}{x_1}+\frac{u}{x_2}+\frac{u}{x_3},\,\,\hat{f}^{\downarrow,{{\theta}_M}}=\frac{p}{x_1}+\frac{p}{x_2}+\frac{p}{x_3},\\
\hat{f}^{\uparrow,{{\eta}_M}}=\frac{r}{x_1}+\frac{r}{x_2}+\frac{r}{x_3},\,\,\hat{f}^{\downarrow,\mathcal{I}_{{\eta}_M}}=\frac{0}{x_1}+\frac{u}{x_2}+\frac{t}{x_3}.$
\end{exa}
\begin{rem}\label{rm42} (i) If $L=[0,1],\,\mathbf{N}=\mathbf{N}_S,{\theta}={\theta}_{M}$ and ${\eta}={\eta}_{M}$, then the inverse $F$-transforms $\hat{f}^{\uparrow,\mathcal{I}_{{\theta}}},\hat{f}^{\downarrow,{\theta}},\hat{f}^{\uparrow,{\eta}}$ and $\hat{f}^{\downarrow,\mathcal{I}_{\eta}}$ become as follows:
\begin{eqnarray*}
{\hat{f}^{\uparrow,\mathcal{I}_{{\theta}_M}}(x)}&=&\bigwedge\limits_{j\in J}\mathcal{I}_{{\theta}_M}(A_j(x), F^{\uparrow,{\theta}_M}_j[f]),\\
{\hat{f}^{\downarrow,{\theta}_M}(x)}&=&\bigvee\limits_{j\in J}(A_j(x)\wedge F^{\downarrow,\mathcal{I}_{{\theta}_M}}_j[f]),\\
{\hat{f}^{\uparrow,{\eta}_M}(x)}&=&\bigwedge\limits_{j\in J}(( 1-A_j(x))\vee F^{\uparrow,\mathcal{I}_{{\eta}_M}}_j[f]),\,\text{and}\\
{\hat{f}^{\downarrow,\mathcal{I}_{{\eta}_M}}(x)}&=&\bigvee\limits_{j\in J}\mathcal{I}_{{\eta}_M}((1-A_j(x)), F^{\downarrow,{\eta}_M}_j[f]),\,\forall\,x\in X,f\in L^X.
\end{eqnarray*}
Obviously $f^{\uparrow,\mathcal{I}_{{\theta}_M}}$ and $f^{\downarrow,{{\theta}_M}}$ coincide with the special cases of inverse upper and lower $F$-transforms proposed in \cite{per,anan,tri}, respectively .\\\\ 
(ii) If $L=[0,1],{\theta}={\theta}_{M}$ and ${\eta}={\eta}_{M}$, then the inverse transforms $\hat{f}^{\uparrow,\mathcal{I}_{{\theta}}},\hat{f}^{\downarrow,{\theta}},\hat{f}^{\uparrow,{\eta}}$ and $\hat{f}^{\downarrow,\mathcal{I}_{\eta}}$ become as follows:
\begin{eqnarray*}
{\hat{f}^{\uparrow,\mathcal{I}_{{\theta}_M}}(x)}&=&\bigwedge\limits_{j\in J}\mathcal{I}_{{\theta}_M}(A_j(x), F^{\uparrow,{\theta}_M}_j[f]),\\
{\hat{f}^{\downarrow,{\theta}_M}(x)}&=&\bigvee\limits_{j\in J}(A_j(x)\wedge F^{\downarrow,\mathcal{I}_{{\theta}_M}}_j[f]),\\
{\hat{f}^{\uparrow,{\eta}_M}(x)}&=&\bigwedge\limits_{j\in J}(\mathbf{N}( A_j(x))\vee F^{\uparrow,\mathcal{I}_{{\eta}_M}}_j[f]),\,\text{and}\\
{\hat{f}^{\downarrow,\mathcal{I}_{{\eta}_M}}(x)}&=&\bigvee\limits_{j\in J}\mathcal{I}_{{\eta}_M}(\mathbf{N}(A_j(x)), F^{\downarrow,{\eta}_M}_j[f]),\,\forall\,x\in X,f\in L^X.
\end{eqnarray*}
Obviously $f^{\uparrow,\mathcal{I}_{{\theta}_M}}$ and $f^{\downarrow,{{\theta}_M}}$ coincide with the special cases of inverse upper and lower $F$-transforms proposed in \cite{per,anan,tri}, respectively.\\\\ 
(iii) If $L=[0,1],{\theta}=\mathcal{T}$ and ${\eta}=\mathcal{S}$, where $\mathcal{T},\mathcal{S}$ are continuous $t$-norm, $t$-conorm with no nontrivial zero divisors, respectively, then the inverse transforms $\hat{f}^{\uparrow,\mathcal{I}_{{\theta}}},\hat{f}^{\downarrow,{\theta}},\hat{f}^{\uparrow,{\eta}}$ and $\hat{f}^{\downarrow,\mathcal{I}_{\eta}}$ become as follows:
\begin{eqnarray*}
{\hat{f}^{\uparrow,\mathcal{I}_{\mathcal{T}}}(x)}&=&\bigwedge\limits_{j\in J}\mathcal{I}_{\mathcal{T}}(A_j(x), F^{\uparrow,\mathcal{T}}_j[f]),\\
{\hat{f}^{\downarrow,\mathcal{T}}(x)}&=&\bigvee\limits_{j\in J}\mathcal{T}(A_j(x), F^{\downarrow,\mathcal{I}_{\mathcal{T}}}_j[f]),\\
{\hat{f}^{\uparrow,\mathcal{S}}(x)}&=&\bigwedge\limits_{j\in J}\mathcal{S}(\mathbf{N}( A_j(x)), F^{\uparrow,\mathcal{I}_{\mathcal{S}}}_j[f]),\,\text{and}\\
{\hat{f}^{\downarrow,\mathcal{I}_{\mathcal{S}}}(x)}&=&\bigvee\limits_{j\in J}\mathcal{I}_{\mathcal{S}}(\mathbf{N}(A_j(x)), F^{\downarrow,\mathcal{S}}_j[f]),\,\forall\,x\in X,f\in L^X.
\end{eqnarray*}
Obviously $f^{\uparrow,\mathcal{I}_{\mathcal{T}}}$ and $f^{\downarrow,{\mathcal{T}}}$ coincide with the of inverse upper and lower $F$-transforms computed with $t$-norm and $R$-implicator proposed in \cite{per,anan,tri}, respectively.
\end{rem} 
From the above, it is clear that some existing inverse $F$-transforms are special cases of the proposed inverse $F$-transforms. Among these, some inverse $F$-transforms coincide with the proposed inverse $F$-transforms and some of the proposed inverse $F$-transforms coincide with the special cases of the existing inverse $F$-transforms. That is to say; the proposed inverse $F$-transforms are more general than some existing ones.\\\\
The following two results are towards the inverse $F$-transforms approximates the original $L$-fuzzy set.
\begin{pro}\label{p41}
{Let $\mathcal{P}$ be an $L$-fuzzy partition of $X$. Then for all $x\in X,f\in L^X$
\begin{itemize}
\item[(i)] $\hat{f}^{\uparrow,\mathcal{I}_{{\theta}}}(x)\geq f(x)$, and
\item[(ii)] $\hat{f}^{\downarrow,{\theta}}(x)\leq f(x)$.
\end{itemize}}
\end{pro}
\textbf{Proof:} (i) Let $x\in X,f\in L^X$. Then from Definition \ref{IT} 
\begin{eqnarray*}
\hat{f}^{\uparrow,\mathcal{I}_{{\theta}}}(x)&=&\bigwedge\limits_{j\in J}\mathcal{I}_{\theta}(A_j(x),F^{\uparrow,{\theta}}_j[f])=\bigwedge\limits_{j\in J}\mathcal{I}_{\theta}(A_j(x),\bigvee\limits_{y\in X}{\theta}(A_j(y),f(y)))\\
&\geq&\bigwedge\limits_{j\in J}\mathcal{I}_{\theta}(A_j(x),{\theta}(A_j(x),f(x)))\geq f(x).
\end{eqnarray*}
Thus $\hat{f}^{\uparrow,\mathcal{I}_{{\theta}}}(x)\geq f(x)$.\\\\
(ii) Let $x\in X$ and $f\in L^X$. Then from Definition {\ref{IT}}
\begin{eqnarray*}
\hat{f}^{\downarrow,{\theta}}(x)&=&\bigvee\limits_{j\in J}{\theta}(A_j(x),F^{\downarrow,\mathcal{I}_{{\theta}}}_j[f])=\bigvee\limits_{j\in J}{\theta}(A_j(x),\bigwedge\limits_{y\in X}\mathcal{I}_{\theta}(A_j(y),f(y)))\\
&\leq&\bigvee\limits_{j\in J}{\theta}(A_j(x),\mathcal{I}_{\theta}(A_j(x),f(x)))\leq f(x).
\end{eqnarray*}
Thus $\hat{f}^{\downarrow,{\theta}}(x)\leq\mathcal{I}_{\theta}(1,f(x))$.
\begin{pro}\label{p42}
{Let $\mathcal{P}$ be an $L$-fuzzy partition of $X$. Then for all $x\in X,f\in L^X$
\begin{itemize}
\item[(i)] $\hat{f}^{\uparrow,{\eta}}(x)\geq f(x)$, and
\item[(ii)] $\hat{f}^{\downarrow,\mathcal{I}_{\eta}}(x)\leq \mathcal{I}_{\eta}(0,f(x))$.
\end{itemize}}
\end{pro}
\textbf{Proof:} (i) Let $x\in X,f\in L^X$. Then from Definition \ref{IT} 
\begin{eqnarray*}
\hat{f}^{\uparrow,{\eta}}(x)&=&\bigwedge\limits_{j\in J}{\eta}(\mathbf{N}(A_j(x)),F^{\uparrow,\mathcal{I}_{\eta}}_j[f])=\bigwedge\limits_{j\in J}{\eta}(\mathbf{N}(A_j(x)),\bigvee\limits_{y\in X}{\mathcal{I}_{\eta}}((\mathbf{N}(A_j(y)),f(y)))\\
&\geq&\bigwedge\limits_{j\in J}{\eta}(\mathbf{N}(A_j(x)),{\mathcal{I}_{\eta}}((\mathbf{N}(A_j(x)),f(x)))\geq f(x).
\end{eqnarray*}
Thus $\hat{f}^{\uparrow,{\eta}}(x)\geq f(x)$.\\\\
(ii) Let $x\in X$ and $f\in L^X$. Then from Definition {\ref{IT}}
\begin{eqnarray*}
\hat{f}^{\downarrow,\mathcal{I}_{\eta}}(x)&=&\bigvee\limits_{j\in J}\mathcal{I}_{\eta}(\mathbf{N}(A_j(x)),F^{\downarrow,{\eta}}_j[f])=\bigvee\limits_{j\in J}\mathcal{I}_{\eta}(\mathbf{N}(A_j(x)),\bigwedge\limits_{y\in X}{\eta}(\mathbf{N}(A_j(y)),f(y)))\\
&\leq&\bigvee\limits_{j\in J}\mathcal{I}_{\eta}(A_j(x),{\eta}(\mathbf{N}(A_j(x)),f(x)))\leq f(x).
\end{eqnarray*}
Thus $\hat{f}^{\downarrow,\mathcal{I}_{\eta}}(x)\leq f(x)$.\\\\
Below, we show that the $L$-fuzzy set $f$ and inverse $F$-transforms have the same $F$-transforms, respectively. Therefore the inverse $F$-transforms of the inverse $F$-transforms is again inverse $F$-transforms, respectively. This can easily follows from the following.
\begin{pro}\label{43}
Let $\mathcal{P}$ be an $L$-fuzzy partition of $X$. Then for all $j\in J,f\in L^X$
\begin{itemize}
\item[(i)] $F_j^{\uparrow,{\theta}}[f]=\bigvee\limits_{x\in X}{\theta}(A_j(x),\hat{f}^{\uparrow,\mathcal{I}_{{\theta}}}(x))$, and
\item[(ii)] $F_j^{\downarrow,\mathcal{I}_{{\theta}}}[f]=\bigwedge\limits_{x\in X}{\mathcal{I}_{{\theta}}}(A_j(x),\hat{f}^{\downarrow,{\theta}}(x))$.
\end{itemize}
\end{pro}
\textbf{Proof:} (i) From Proposition \ref{p41}(i), $\hat{f}^{\uparrow,\mathcal{I}_{{\theta}}}(x)\geq f(x),\,\forall\,x\in X$. It follows from Definition \ref{UFT} that
\begin{eqnarray*}
F_j^{\uparrow,{\theta}}[f]=\bigvee\limits_{x\in X}{\theta}(A_j(x),{f}(x))&\leq&\bigvee\limits_{x\in X}{\theta}(A_j(x),\hat{f}^{\uparrow,\mathcal{I}_{{\theta}}}(x))~\text{and}\\
{\theta}(A_j(x),\hat{f}^{\uparrow,\mathcal{I}_{{\theta}}}(x))&=&{\theta}(A_j(x),\bigwedge\limits_{k\in J}\mathcal{I}_{\theta}(A_k(x),F^{\uparrow,{\theta}}_k[f])\\
 &\leq&{\theta}(A_j(x),\mathcal{I}_{\theta}(A_j(x),F^{\uparrow,{\theta}}_j[f]) \\
 &\leq& {F^\uparrow_j[f].}
 \end{eqnarray*}
 Thus $\bigvee\limits_{x\in X}{\theta}(A_j(x),\hat{f}^{\uparrow,\mathcal{I}_{{\theta}}}(x))\leq F^{\uparrow,{\theta}}_j[f]$ or $F_j^{\uparrow,{\theta}}[f]=\bigvee\limits_{x\in X}{\theta}(A_j(x),\hat{f}^{\uparrow,\mathcal{I}_{{\theta}}}(x))$.\\\\
 (ii) From Proposition \ref{p41}(ii), $\hat{f}^\downarrow(x)\leq f(x),\,\forall\,x\in X$. It follows from Definition \ref{UFT} that
\begin{eqnarray*}
F_j^{\downarrow,\mathcal{I}_{{\theta}}}[f]=\bigwedge\limits_{x\in X}\mathcal{I}_{\theta}(A_j(x),{f}(x))&\geq&\bigwedge\limits_{x\in X}\mathcal{I}_{\theta}(A_j(x),\hat{f}^{\downarrow,{\theta}}(x))~\text{and}\\
 \mathcal{I}_{\theta}(A_j(x),\hat{f}^{\downarrow,{\theta}}(x))&=&\mathcal{I}_{\theta}(A_j(x),\bigvee\limits_{k\in J}{\theta}(A_k(x),F^{\downarrow,\mathcal{I}_{{\theta}}}_k[f]))\\
 &\geq&\mathcal{I}_{\theta}(A_j(x),{\theta}(A_j(x),F^{\downarrow,\mathcal{I}_{{\theta}}}_j[f]))\\
 &\geq&F^{\downarrow,\mathcal{I}_{{\theta}}}_j[f].
 \end{eqnarray*}
 Thus $\bigwedge\limits_{x\in X}\mathcal{I}_{\theta}(A_j(x),\hat{f}^{\downarrow,{\theta}}(x))\geq F^{\downarrow,\mathcal{I}_{{\theta}}}_j[f]$ or $F_j^{\downarrow,\mathcal{I}_{{\theta}}}[f]=\bigwedge\limits_{x\in X}\mathcal{I}(A_j(x),\hat{f}^{\downarrow,{\theta}}(x)).$
\begin{pro}\label{44}
Let $\mathcal{P}$ be an $L$-fuzzy partition of $X$. Then for all $j\in J,f\in L^X$
\begin{itemize}
\item[(i)] $F_j^{\uparrow,\mathcal{I}_{{\eta}}}[f]=\bigvee\limits_{x\in X}{\mathcal{I}_{{\eta}}}(\mathbf{N}(A_j(x)),\hat{f}^{\uparrow,{\eta}}(x))$, and
\item[(ii)] $F_j^{\downarrow,{\eta}}[f]=\bigwedge\limits_{x\in X}{\eta}(\mathbf{N}(A_j(x)),\hat{f}^{\downarrow,\mathcal{I}_{\eta}}(x))$.
\end{itemize}
\end{pro}
\textbf{Proof:} (i) From Proposition \ref{p42}(i), $\hat{f}^{\uparrow,{\eta}}(x)\geq f(x),\,\forall\,x\in X$. It follows from Definition \ref{UFT} that
\begin{eqnarray*}
F_j^{\uparrow,\mathcal{I}_{\eta}}[f]=\bigvee\limits_{x\in X}\mathcal{I}_{\eta}(\mathbf{N}(A_j(x)),{f}(x))&\leq&\bigvee\limits_{x\in X}\mathcal{I}_{\eta}(\mathbf{N}(A_j(x)),\hat{f}^{\uparrow,{\eta}}(x))~\text{and}\\
\mathcal{I}_{\eta}(\mathbf{N}(A_j(x)),\hat{f}^{\uparrow,{\eta}}(x))&=&\mathcal{I}_{\eta}(\mathbf{N}(A_j(x)),\bigwedge\limits_{k\in J}{\eta}(\mathbf{N}(A_k(x)),F^{\uparrow,\mathcal{I}_{\eta}}_k[f]))\\
 &\leq&\mathcal{I}_{\eta}(\mathbf{N}(A_j(x)),{\eta}(\mathbf{N}(A_j(x)),F^{\uparrow,\mathcal{I}_{\eta}}_j[f])) \\
 &\leq&{F^{\uparrow,\mathcal{I}_{\eta}}_j[f]}.
 \end{eqnarray*}
 Thus $\bigvee\limits_{x\in X}\mathcal{I}_{\eta}(\mathbf{N}(A_j(x)),\hat{f}^{\uparrow,{\eta}}(x))\leq F^{\uparrow,\mathcal{I}_{\eta}}_j[f]$ or $F_j^{\uparrow,\mathcal{I}_{\eta}}[f]=\bigvee\limits_{x\in X}\mathcal{I}_{\eta}(A_j(x),\hat{f}^{\uparrow,{\eta}}(x))$.\\\\
 (ii) From Proposition \ref{p42}(ii), $\hat{f}^{\downarrow,\mathcal{I}_{\eta}}(x)\leq f(x),\,\forall\,x\in X$. It follows from Definition \ref{UFT} that
\begin{eqnarray*}
F_j^{\downarrow,{\eta}}[f]=\bigwedge\limits_{x\in X}{\eta}(\mathbf{N}(A_j(x)),{f}(x))&\geq&\bigwedge\limits_{x\in X}{\eta}(\mathbf{N}(A_j(x)),\hat{f}^{\downarrow,{\mathcal{I}_{\eta}}}(x))~\text{and}\\
{\eta}(\mathbf{N}(A_j(x)),\hat{f}^{\downarrow,\mathcal{I}_{\eta}}(x))&=&{\eta}(\mathbf{N}(A_j(x)),\bigvee\limits_{k\in J}\mathcal{I}_{\eta}(\mathbf{N}(A_k(x)),F^{\downarrow,{\eta}}_k[f]))\\
 &\geq&{\eta}(\mathbf{N}(A_j(x)),\mathcal{I}_{\eta}(\mathbf{N}(A_j(x)),F^{\downarrow,{\eta}}_j[f]))\\
 &\geq& F^{\downarrow,{\eta}}_j[f].
 \end{eqnarray*}
 Thus $\bigwedge\limits_{x\in X}{\eta}(\mathbf{N}(A_j(x)),\hat{f}^{\downarrow,\mathcal{I}_{\eta}}(x))\geq F^{\downarrow,{\eta}}_j[f]$ or $F_j^{\downarrow,{\eta}}[f]=\bigwedge\limits_{x\in X}\mathcal{I}(A_j(x),\hat{f}^{\downarrow,\mathcal{I}_{\eta}}(x)).$
\section{Axiomatic approaches of $F$-transforms}
 In \cite{jiri}, the axiomatic approaches of the direct $F$-transforms computed with $t$-norm and $R$-implicator were studied in detail. This section focuses on the axiomatic characterizations of the direct $F$-transforms computed with respect to ${\theta},{\eta},\mathcal{I}_{\eta},\mathcal{I}_{{\theta}}$, respectively by some independent axioms. Also, we first present the axioms for each direct $F$-transform that guarantee the existence of an $L$-fuzzy partition that produces the same $F$-transform. Now, we begin with the following.\\\\
 For any $f\in L^X$ and for an {$L$-fuzzy} partition $\mathcal{P}$, it can be seen that the direct $F^{\uparrow,{\theta}},F^{\downarrow,{\eta}},F^{\uparrow,\mathcal{I}_{\eta}}$ and $F^{\uparrow,\mathcal{I}_{{\theta}}}$-transforms {induce} the maps $F^{\uparrow,{\theta}}_{ \mathcal{P}},F^{\downarrow,{\eta}}_{ \mathcal{P}},F^{\uparrow,\mathcal{I}_{\eta}}_{ \mathcal{P}},F^{\uparrow,\mathcal{I}_{{\theta}}}_{ \mathcal{P}}:L^X\rightarrow L^J$ such that
 \begin{eqnarray*}
 F^{\uparrow,{\theta}}_{\mathcal{P}}[f](j) &=& F^{\uparrow,{\theta}}_{j}[f],\,\, F^{\downarrow,{\eta}}_{\mathcal{P}}[f](j) = F^{\downarrow,{\eta}}_{j}[f],\\
 F^{\uparrow,\mathcal{I}_{\eta}}_{\mathcal{P}}[f](j) &=& F^{\uparrow,\mathcal{I}_{\eta}}_{j}[f],\,\, F^{\downarrow,\mathcal{I}_{{\theta}}}_{\mathcal{P}}[f](j) = F^{\downarrow,\mathcal{I}_{{\theta}}}_{j}[f],\,\text{respectively}.
 \end{eqnarray*}
 Now, we introduce the concepts of $L$-fuzzy upper and lower transformation systems with respect to overlap and grouping maps $\theta$ and $\eta$ ( a co-residual and residual implicators $\mathcal{I}_\eta$ and $\mathcal{I}_\theta$ induced by grouping and overlap maps ${\eta}$ and $\theta$), respectively.
 \begin{def1}\label{UTS} Let $X$ be a nonempty set, ${{\theta}}$ be an overlap map and $\mathcal{I}_{\eta}$ be a co-residual implicator over $L$. Then a system $\mathcal{U}_F = (X,Y, u, U_F)$, where $F={\theta}\,\,\text{or}\,\,\mathcal{I}_{\eta}$ and
 \begin{itemize}
\item[1.] $Y$ is a nonempty set,
\item[2.] $u : X \rightarrow Y $ is an onto map,
\item[3.] $U_F : L^X \rightarrow L^Y$ is a map, where
 \begin{itemize}
\item[(i)] for all $\lbrace f_k : k \in J\rbrace$ $\subseteq L^X$,
$U_F[\bigvee\limits_{k\in J}f_k]= \bigvee\limits_{k\in J}U_F[f_k]$,
\item[(ii)] for all $\textbf{u},f \in L^X$, $U_F[F(\textbf{u},f)]= {F}(\textbf{u}, U_F[f])$,
\item[(iii)] for all $x\in X,\,y\in Y$, $U_F[1_{\lbrace x\rbrace}](y) = 1$ iff $y = u(x)$,
\end{itemize}
 \end{itemize}
is called an {\bf $L$-fuzzy upper transformation system} on $X$ with respect to $F$.
\end{def1}
\begin{def1}\label{LTS} {Let $X$ be a nonempty set, ${\eta},\mathcal{I}_\theta$ and $\mathbf{N}$ be grouping map, residual implicator and negator over $L$, respectively. Then
system $\mathcal{H}_{{F}} = (X, Y, v,H_F)$, where $F={\eta}\,\,\text{or}\,\,\mathcal{I}_{{\theta}}$ and
\begin{itemize}
\item[1.] $Y$ is a nonemptyy set,
\item[2.] $v : X \rightarrow Y $ is an onto map,
\item[3.]\label{prop} $H_F : L^X \rightarrow L^Y$ is a map, where
\begin{itemize}
\item[(i)] for all $\lbrace f_ k: k\in J\rbrace$ $\subseteq L^X$, $H_F[\bigwedge\limits_{k\in J}f_k](y) = \bigwedge\limits_{k\in J}H_F[f_k](y),$
\item[(ii)]{for all $\textbf{u},f \in L^X$, ${H_F[F(\textbf{u}, f)]}= F(\textbf{u}, H_F[f])$,} and
\item[(iii)] for $y\in Y$ and $x\in X$, ${(\mathbf{N}( H_F[{\mathbf{N}}(1_{\lbrace x\rbrace})]))(y)} = 1$ iff $y = v(x)$,
\end{itemize}
\end{itemize}
is called an {\bf $L$-fuzzy lower transformation system} on $X$ with respect to ${F}$}.
\end{def1}
The $L$-fuzzy upper transformation system with respect to a $t$-norm and the $L$-fuzzy lower transformation system with respect to an $R$-implicator proposed in \cite{jiri,tri} are special cases of $\mathcal{U}_\theta$ and $\mathcal{H}_{\mathcal{I}_\theta}$, respectively. Also, the $L$-fuzzy lower transformation system with respect to an $S$-implicator proposed in \cite{tri} is a special case of $\mathcal{H}_{\eta}$. The $L$-fuzzy upper transformation system with respect to $\mathcal{U}_{\mathcal{I}_{\eta}}$ is a new definition.
\begin{exa}\label{exa51}
Let $X$ be a nonempty set and $id:X\rightarrow X$ be an identity map. Now, we define maps $U_F,H_{F'}:L^X\rightarrow L^X$ such that
$U_{F}[f](x)=f(x),H_{F'}[f](x)=f(x), x\in X$, where $F={\theta}\,or\,\mathcal{I}_{{\eta}}$ and $F'={\eta}\,or\,\mathcal{I}_{{\theta}}$. Then for all $\{{f_k:k\in J}\}\subseteq L^X$,
$U_{F}[\bigvee\limits_{k\in J}f_k]=\bigvee\limits_{k\in J}f_k=\bigvee\limits_{k\in J}U_F[f_k]$ and
$H_{F'}[\bigwedge\limits_{k\in J}f_k]=\bigwedge\limits_{k\in J}f_k=\bigwedge\limits_{k\in J}H_{F'}[f_k].$
Now, let $\textbf{u},f\in L^X$. Then
$U_F[F(\textbf{u}, f)]={F}(\textbf{u},U_F[f])$ and
$H_{F'}[F'(\textbf{u}, f)]={F'}(\textbf{u},H_{F'}[f])$. Finally, let $x,z\in X$. Then
 $U_F[1_{\lbrace x\rbrace}](z)=U_F[1_{\lbrace x\rbrace}](z)= 1_{\lbrace x\rbrace}(z)=1$ iff $x=z$ and
 $(\mathbf{N}(H_{F'}[\mathbf{N}(1_{\lbrace x\rbrace})]))(z)=\mathbf{N}(H_{F'}[\mathbf{N}(1_{\lbrace x\rbrace})](z))= \mathbf{N}(\mathbf{N}(1_{\lbrace x\rbrace})(z))=1$ iff $x=z.$
Thus $U_F(1_{\lbrace x\rbrace})(z)=1,(\mathbf{N}(H_{F'}[\mathbf{N}(1_{\lbrace x\rbrace})]))(z)=1$ iff $z=id(x)$. Hence $\mathcal{U}_{{F}}=(X,X, id, U_F)$ and $\mathcal{H}_{{F'}}=(X,X, id, H_{F'})$ are $L$-fuzzy upper and lower transformation systems on $X$ with respect to $F$ and $F'$, respectively.
\end{exa}
\begin{rem}\label{rm52} (i) If $L=[0,1],\,\mathbf{N}=\mathbf{N}_S,{\theta}={\theta}_{M},{\eta}={\eta}_{M},\mathcal{I}_{\eta}=\mathcal{I}_{{\eta}_M}$ and $\mathcal{I}_{{\theta}}=\mathcal{I}_{{\theta}_M}$. Then $\mathcal{U}_{\theta_M}$ and $\mathcal{H}_{\mathcal{I}_{{\theta}_M}}$ coincide with the special cases of the $L$-fuzzy upper and lower transformation systems proposed in \cite{jiri,tri}, respectively. Also, $\mathcal{H}_{{\eta}_M}$ coincides with the special case of the $L$-fuzzy lower transformation system proposed in \cite{tri}.\\\\ 
(ii) If $L=[0,1],{\theta}={\theta}_{M},{\eta}={\eta}_{M},\mathcal{I}_{\eta}=\mathcal{I}_{{\eta}_M}$ and $\mathcal{I}_{\eta}=\mathcal{I}_{{\theta}_M}$. Then $\mathcal{U}_{\theta_M}$ and $\mathcal{H}_{\mathcal{I}_{{\theta}_M}}$ coincide with the special cases of the $L$-fuzzy upper and lower transformation systems proposed in \cite{jiri,tri}, respectively. Also, $\mathcal{H}_{{\eta}_M}$ coincides with the special case of the $L$-fuzzy lower transformation system proposed in \cite{tri}.\\\\
(iii) If $L=[0,1],{\theta}=\mathcal{T}$ and ${\eta}=\mathcal{S}$, where $\mathcal{T},\mathcal{S}$ are continuous $t$-norm, $t$-conorm with no nontrivial zero divisors, respectively. Then $\mathcal{U}_{\mathcal{T}}$ and $\mathcal{H}_{\mathcal{I}_{\mathcal{T}}}$ coincide with the $L$-fuzzy upper and lower transformation systems with respect to $t$-norm and $R$-implicator proposed in \cite{jiri,tri}, respectively. Also, $\mathcal{H}_{\mathcal{S}}$ coincides with the $L$-fuzzy lower transformation system with respect to $S$-implicator proposed in \cite{tri}.
\end{rem} 
From the above remark, it is clear that some existing $L$-fuzzy transformation systems are special cases of the proposed $L$-fuzzy transformation systems. Among these, some $L$-fuzzy transformation systems coincide with the proposed $L$-fuzzy transformation systems, and some proposed $L$-fuzzy transformation systems coincide with the special cases of the existing $L$-fuzzy transformation systems. That is to say; the proposed $L$-fuzzy transformation systems are more extended form than some existing ones.\\\\
The following shows a close connection of the {$L$-fuzzy} transformation systems with the $F$-transforms. To do this, we need some results, which are given by the following proposition.
\begin{pro}
Let $\mathbf{N}$ be a negator, ${\theta},{\eta}$ be overlap and grouping maps with neutral elements $0,1$, respectively. In addition, let $\mathbf{N}_{\mathcal{I}_{\eta}},\mathbf{N}_{\mathcal{I}_{\eta}}$ be involutive negators. Then for all $f\in L^X$
\begin{itemize}
\item[(i)] $f=\bigvee\limits_{x\in X}{\theta}(\textbf{f(x)},1_{\{x\}}),f=\bigwedge\limits_{x\in X}{\eta}(\textbf{f(x)},\mathbf{N}(1_{\{x\}}))$, and
\item[(ii)] $f=\bigvee\limits_{x\in X}\mathcal{I}_{\eta}(\mathbf{N}_{\mathcal{I}_{\eta}}(\textbf{f(x)}),1_{\{x\}}),f=\bigwedge\limits_{x\in X}\mathcal{I}_{\theta}(\mathbf{N}_{\mathcal{I}_{{\theta}}}(\textbf{f(x)}),\mathbf{N}_{\mathcal{I}_{{\theta}}}(1_{\{x\}}))$.
\end{itemize}
\end{pro}
\textbf{Proof:} (i) Let $y\in X,f\in L^X$. Then
\begin{eqnarray*}
f(y)&=&\bigvee\limits_{x\in X}{\theta}(\textbf{f(x)},1_{\{x\}})(y)=\bigvee\limits_{x\in X}{\theta}({f(x)},1_{\{x\}}(y))\\
&=&{\theta}({f(x)},1_{\{y\}}(y))\vee \bigvee\limits_{x\neq y\in X}{\theta}({f(y)},1_{\{x\}}(y))\\
&=&{\theta}({f(y)},1)={f(y)}.
\end{eqnarray*}
Thus $f=\bigvee\limits_{x\in X}{\theta}(\textbf{f(x)},1_{\{x\}})$ and
\begin{eqnarray*}
f(y)&=&\bigwedge\limits_{x\in X}{\eta}(\textbf{f(x)},\mathbf{N}(1_{\{x\}}))(y)= \bigwedge\limits_{x\in X}{\eta}({f(x)},\mathbf{N}(1_{\{x\}})(y))\\
&=&{\eta}({f(y)},\mathbf{N}(1_{\{y\}}(y)))\vee \bigwedge\limits_{x\neq y\in X}{\eta}({f(y)},\mathbf{N}(1_{\{x\}}(y)))\\
&=&{\eta}({f(y)},0)= f(y).
\end{eqnarray*}
Thus $f=\bigwedge\limits_{x\in X}{\eta}(\textbf{f(x)},\mathbf{N}(1_{\{x\}}))$.\\\\
(ii) Let $y\in X,f\in L^X$. Then
\begin{eqnarray*}
f(y)&=&\bigvee\limits_{x\in X}\mathcal{I}_{\eta}(\mathbf{N}_{\mathcal{I}_{\eta}}(\textbf{f(x)}),1_{\{x\}})(y)=\bigvee\limits_{x\in X}\mathcal{I}_{\eta}(\mathbf{N}_{\mathcal{I}_{\eta}}({f(x)}),1_{\{x\}}(y))\\
&=&\mathcal{I}_{\eta}(\mathbf{N}_{\mathcal{I}_{\eta}}({f(x)}),1_{\{y\}}(y))\vee \bigvee\limits_{x\neq y\in X}\mathcal{I}_{\eta}(\mathbf{N}_{\mathcal{I}_{\eta}}({f(x)}),1_{\{x\}}(y))\\
&=&\mathcal{I}_{\eta}(\mathbf{N}_{\mathcal{I}_{\eta}}({f(y)}),1)=\mathbf{N}_{\mathcal{I}_{\eta}}(\mathbf{N}_{\mathcal{I}_{\eta}}({f(y)})= f(y).
\end{eqnarray*}
Thus $f=\bigvee\limits_{x\in X}\mathcal{I}_{\eta}(\mathbf{N}_{\mathcal{I}_{\eta}}(\textbf{f(x)}),1_{\{x\}})$ and 
\begin{eqnarray*}
f(y)&=&\bigwedge\limits_{x\in X}\mathcal{I}_{\theta}(\mathbf{N}_{\mathcal{I}_{{\theta}}}(\textbf{f(x)}),\mathbf{N}_{\mathcal{I}_{{\theta}}}(1_{\{x\}}))(y)=\bigwedge\limits_{x\in X}\mathcal{I}_{\theta}(\mathbf{N}_{\mathcal{I}_{{\theta}}}({f(x)}),\mathbf{N}_{\mathcal{I}_{{\theta}}}(1_{\{x\}}(y)))\\
&=&\mathcal{I}_{\theta}(\mathbf{N}_{\mathcal{I}_{{\theta}}}({f(y)}),\mathbf{N}_{\mathcal{I}_{{\theta}}}(1_{\{y\}}(y)))\wedge\bigwedge\limits_{x\neq y\in X}\mathcal{I}_{\theta}(\mathbf{N}_{\mathcal{I}_{{\theta}}}({f(x)}),\mathbf{N}_{\mathcal{I}_{{\theta}}}(1_{\{x\}}(y)))\\
&=&\mathcal{I}_{\theta}(\mathbf{N}_{\mathcal{I}_{{\theta}}}({f(y)}),0)=\mathbf{N}_{\mathcal{I}_{{\theta}}}(\mathbf{N}_{\mathcal{I}_{{\theta}}}({f(y)}))= f(y).
\end{eqnarray*}
Now, we have the following.
\begin{pro}\label{p51} Let ${{\theta}}$ be an overlap map $L$. Then the following statements are equivalent:
\begin{itemize}
\item[(i)] $\mathcal{U}_{{\theta}}=(X,Y, u, U_{{\theta}})$ is an $L$-fuzzy upper transformation system on $X$ determined by an overlap map ${{\theta}}$ and {$Y\subseteq X$}.
\item[(ii)] There exists an $L$-fuzzy partition $\mathcal{P}$ of $X$ {indexed by $Y$} such that $u(x) = y$ iff $x\in core( A_y)$ and $U_{{\theta}} = F^{\uparrow,{\theta}}_{ \mathcal{P}}$.
 \end{itemize}
\end{pro}
\textbf{Proof:} Let ${\mathcal{U}_{{\theta}}=(X,Y,u,U_{{\theta}})}$ be an {$L$-fuzzy} upper transformation system on $X$ {determined by ${{\theta}}$}. Also, let $\mathcal{P} = \lbrace A_y : y\in Y\rbrace$ such that for all $y\in Y$, $A_y \in L^X$ is given by $A_y(x) = U_{{\theta}}[ 1_{\lbrace x\rbrace}](y)$, $x\in X$. Now, from Definition \ref{UTS}(iii), $A_{u(x)}(x) = U_{{\theta}}[ 1_{\lbrace x\rbrace}](u(x)) = 1$, or that, $x \in core( A_{u(x)})$. Further, for $y, z \in Y,t \in core ( A_y) \cap core (A_z),U_{{\theta}}[ 1_{\lbrace t\rbrace}](y)= 1= U_{{\theta}}[ 1_{\lbrace t\rbrace}](z)$, i.e., $ A_y(t) = 1 = A_z(t)$ iff $y = u(t) = z$. Thus $\lbrace core(A_y):y\in Y\rbrace$ is a partition of $X$ and therefore $\mathcal{P}$ is an {$L$-fuzzy} partition of $X$. Now, for all $y\in Y$ and $f\in L^X$
\begin{eqnarray*}
F^{\uparrow,{\theta}} _{\mathcal{P}}[f](y)&=&\bigvee\limits_{x\in X}{\theta}( A_y(x),f(x))\\
&=&\bigvee\limits_{x\in X}{\theta}( U_{{\theta}}[1_{\lbrace x\rbrace}](y), f(x))\\
&=&\bigvee\limits_{x\in X}{\theta}( f(x),U_{{\theta}}[1_{\lbrace x\rbrace}](y))\\
&=&{\bigvee\limits_{x\in X}U_{{\theta}}[{\theta}({\textbf{f(x)}}, 1_{\lbrace x\rbrace})](y)}\\
&=&{U_{{\theta}}[\bigvee\limits_{x\in X}{\theta}({\textbf{f(x)}}, 1_{\lbrace x\rbrace})](y)}\\
&=&U_{{\theta}}[f](y).
\end{eqnarray*}
Thus $U_{{\theta}}=F^{\uparrow,{\theta}}_{\mathcal{P}}$. Conversely, let $\mathcal{P}=\{A_y\in L^X: y\in Y\}$ be an {$L$-fuzzy} partition of base set $X\neq\emptyset$. Let us define a map $u: X \rightarrow Y$ such that $u(x) = y$ iff $x \in core(A_{y})$. Further, let ${{\theta}}$ be an overlap map with neutral element $1$ and {$U_{{\theta}}=F^{\uparrow,{\theta}}_{\mathcal{P}}$}. Then for all $y\in Y,x\in X$
$U_{{\theta}}[ 1_{\lbrace x\rbrace}](y)= F^{\uparrow,{\theta}}_{\mathcal{P}}[ 1_{\lbrace x\rbrace}](y)= \bigvee_{z\in X}{\theta}(A_{y}(z),\ 1_{\lbrace x\rbrace}(z))={\theta}(A_{y}(x),1))=A_{y}(x).$ Thus $ U_{{\theta}}[1_{\lbrace x\rbrace}](y)=1$ iff $A_{y}(x)=1$ iff $v(x)=y$. From Propositions \ref{p35} and \ref{p36}, $(X,Y,u,U_{{\theta}})$ is an {$L$-fuzzy} upper transformation system on $X$ {determined by ${\theta}$}.
\begin{pro}\label{p52} Let $\mathcal{I}_{\eta}$ be an $EP$-co-residual implicator over $L$ such that $\mathbf{N}_{\mathcal{I}_{\eta}}$ is an involutive negator. Then the following statements are equivalent:
\begin{itemize}
\item[(i)] $\mathcal{U}_{\mathcal{I}_{{\eta}}}=(X,Y, u, U_{\mathcal{I}_{{\eta}}})$ is an $L$-fuzzy upper transformation system on $X$ determined by a co-residual implicator $\mathcal{I}_{\eta}$ and {$Y\subseteq X$}.
\item[(ii)] There exists an $L$-fuzzy partition $\mathcal{P}$ of $X$ {indexed by $Y$} such that $u(x) = y$ iff $x\in core( A_y)$ and $U_{\mathcal{I}_{{\eta}}} = F^{\uparrow,\mathcal{I}_{\eta}}_{ \mathcal{P}}$.
\end{itemize}
\end{pro}
\textbf{Proof:} Let ${\mathcal{U}_{\mathcal{I}_{\eta}}=(X,Y,u,U_{\mathcal{I}_{{\eta}}})}$ be an {$L$-fuzzy} upper transformation system on $X$ {determined by $\mathcal{I}_{\eta}$}. Also, let $\mathcal{P} = \lbrace A_y : y\in Y\rbrace$ such that for all $y\in Y$, $A_y \in L^X$ is given by {$A_y(x) = U_{\mathcal{I}_{{\eta}}}[1_{\lbrace x\rbrace}](y)$}, $x\in X$. Now, from Definition \ref{UTS}(iii), $A_{u(x)}(x) = U_{\mathcal{I}_{{\eta}}}[ 1_{\lbrace x\rbrace}](u(x))= 1$, or that, $x \in core( A_{u(x)})$. Further, {for $y, z \in Y,t \in core ( A_y) \cap core (A_z)$}, $ U_{\mathcal{I}_{{\eta}}}[1_{\lbrace t\rbrace}](y)= 1 = U_{\mathcal{I}_{{\eta}}}[ 1_{\lbrace t\rbrace}](z)$, i.e., $ A_y(t) = 1 = A_z(t)$ iff $y = u(t) = z$. Thus $\lbrace core(A_y):y\in Y\rbrace$ is a partition of $X$ and therefore $\mathcal{P}$ is an {$L$-fuzzy} partition of $X$. Now, for all $y\in Y$ and $f\in L^X$
\begin{eqnarray*}
F^{\uparrow,\mathcal{I}_{\eta}} _{\mathcal{P}}[f](y)&=&\bigvee\limits_{x\in X}\mathcal{I}_{\eta}( \mathbf{N}_{\mathcal{I}_{\eta}}(A_y(x)),f(x))\\
&=&\bigvee\limits_{x\in X}\mathcal{I}_{\eta}( \mathbf{N}_{\mathcal{I}_{\eta}}(A_y(x)),\mathbf{N}_{\mathcal{I}_{\eta}}(\mathbf{N}_{\mathcal{I}_{\eta}}(f(x))))\\
&=&\bigvee\limits_{x\in X}\mathcal{I}_{\eta}( \mathbf{N}_{\mathcal{I}_{\eta}}(f(x)),A_y(x))\\
&=&\bigvee\limits_{x\in X}\mathcal{I}_{\eta}( \mathbf{N}_{\mathcal{I}_{\eta}}(f(x)),U_{\mathcal{I}_{{\eta}}}[1_{\lbrace x\rbrace}](y))\\
&=&\bigvee\limits_{x\in X}U_{\mathcal{I}_{{\eta}}}[\mathcal{I}_{\eta}( \mathbf{N}_{\mathcal{I}_{\eta}}(\textbf{f(x)}),1_{\lbrace x\rbrace})](y)\\
&=&U_{\mathcal{I}_{{\eta}}}[\bigvee\limits_{x\in X}\mathcal{I}_{\eta}( \mathbf{N}_{\mathcal{I}_{\eta}}(\textbf{f(x)}),1_{\lbrace x\rbrace})](y)\\
&=&U_{\mathcal{I}_{{\eta}}}[f](y).
\end{eqnarray*}
Thus $U_{\mathcal{I}_{{\eta}}}=F^{\uparrow,\mathcal{I}_{\eta}}_{\mathcal{P}}$. Conversely, let $\mathcal{P}=\{A_y\in L^X: y\in Y\}$ be an {$L$-fuzzy} partition of base set $X\neq\emptyset$. Let us define a map $u: X \rightarrow Y$ such that $u(x) = y$ iff $x \in core(A_{y})$. Further, let $\mathcal{I}_{\eta}$ be a co-residual implicator such that ${\mathbf{N}_{\mathcal{I}_{\eta}}}(\cdot)=\mathcal{I}_{\eta}(\cdot,1)$ is an involutive negator) and {$U_{\mathcal{I}_{{\eta}}}=F^{\uparrow,\mathcal{I}_{\eta}}_{\mathcal{P}}$}. Then for all $y\in Y,x\in X,\,U_{\mathcal{I}_{{\eta}}}[1_{\lbrace x\rbrace}](y)= F^{\uparrow,\mathcal{I}_{\eta}}_{\mathcal{P}}[ 1_{\lbrace x\rbrace}](y)= \bigvee_{z\in X}\mathcal{I}_{\eta}(\mathbf{N}_{\mathcal{I}_{\eta}}(A_{y}(z)), 1_{\lbrace x\rbrace}))(z))= \bigwedge_{z\in X}\mathcal{I}_{\eta}(\mathbf{N}_{\mathcal{I}_{\eta}}(A_{y}(z)), 1_{\lbrace x\rbrace}(z))\linebreak=\mathcal{I}_{\eta}(\mathbf{N}_{\mathcal{I}_{\eta}}(A_{y}(x)),1))=\mathbf{N}_{\mathcal{I}_{\eta}}(\mathbf{N}_{\mathcal{I}_{\eta}}(A_y(x)))=A_{y}(x).$ Thus $U_{\mathcal{I}_{{\eta}}}[ 1_{\lbrace x\rbrace}](y)=1$ iff $A_{y}(x)=1$ iff $u(x)=y$. From Propositions \ref{p35} and \ref{p36}, $(X,Y,u,U_{\mathcal{I}_{{\eta}}})$ is an {$L$-fuzzy} upper transformation system on $X$ {determined by ${\mathcal{I}_{\eta}}$}.
\begin{pro}\label{p53}
Let ${\eta}$ be an $EP$-grouping map with neutral element $0$ over $L$ such that $\mathbf{N}$ be an involutive negator. Then the following statements are equivalent:
\begin{itemize}
\item[(i)] $\mathcal{H}_{{\eta}}=(X,Y, v, H_{\eta})$ is an $L$-fuzzy lower transformation system on $X$ determined by ${\eta}$ and {$Y\subseteq X$}.
\item[(ii)] There exists an $L$-fuzzy partition $\mathcal{P}$ of $X$ {indexed by $Y$}, such that $v(x) = y$ iff $x\in core(A_y)$ and $H_{\eta} = F^{\downarrow,{\eta}}_{\mathcal{P}}$.
\end{itemize}
\end{pro}
\textbf{Proof:} Let ${\mathcal{H}_{{\eta}}=(X,Y,v,H_{\eta})}$ be an {$L$-fuzzy} lower transformation system on $X$ {determined by $\eta$}. Also, let $\mathcal{P} = \lbrace A_y : y\in Y\rbrace$ such that for all $y\in Y$, $A_y \in L^X$ is given by $A_y(x) = {\mathbf{N}}(H_{\eta}[\mathbf{N}( 1_{\lbrace x\rbrace})])(y)$, $x\in X$. Now, from Definition \ref{LTS}(iii), $A_{v(x)}(x) = {({\mathbf{N}}( H_{\eta}[{\mathbf{N}}( 1_{\lbrace x\rbrace})]))(v(x))} = 1$, or that, $x \in core( A_{v(x)})$. Further, for $y, z \in Y,t \in core ( A_y) \cap core (A_z),({\mathbf{N}}( H_{\eta}[{\mathbf{N}}( 1_{\lbrace t\rbrace})]))(y)= 1
=({\mathbf{N}}( H_{\eta}[{\mathbf{N}}( 1_{\lbrace t\rbrace})]))(z)$, i.e., $ A_y(t) = 1 = A_z(t)$ iff $y = v(t) = z$. Thus $\lbrace core(A_y):y\in Y\rbrace$ is a partition of $X$ and therefore $\mathcal{P}$ is an {$L$-fuzzy} partition of $X$. Now, for all $y\in Y$ and $f\in L^X$
\begin{eqnarray*}
F^{\downarrow,{\eta}} _{\mathcal{P}}[f](y)&=&\bigwedge\limits_{x\in X}{\eta}( \mathbf{N}(A_y(x)),f(x))\\
&=&\bigwedge\limits_{x\in X}{\eta}( H_{\eta}[{\mathbf{N}}( 1_{\lbrace x\rbrace})](y), f(x))\\
&=&\bigwedge\limits_{x\in X}{\eta}(f(x),H_{\eta}[{\mathbf{N}}( 1_{\lbrace x\rbrace})](y))\\
&=&{\bigwedge\limits_{x\in X}H_{\eta}[{\eta}({\textbf{f(x)}}, {\mathbf{N}}( 1_{\lbrace x\rbrace}))](y)}\\
&=&{H_{\eta}[\bigwedge\limits_{x\in X}{\eta}({\textbf{f(x)}}, {\mathbf{N}}( 1_{\lbrace x\rbrace}))](y)}\\
&=&H_{\eta}[f](y).
\end{eqnarray*}
Thus $H_{\eta}=F^{\downarrow,{\eta}}_{\mathcal{P}}$. Conversely, let $\mathcal{P}=\{A_y\in L^X: y\in Y\}$ be an {$L$-fuzzy} partition of base set $X\neq\emptyset$. Let us define a map $v: X \rightarrow Y$ such that $v(x) = y$ iff $x \in core(A_{y})$. Further, let ${\eta}$ be a grouping map with neutral element $0$, ${\mathbf{N}}$ be an involutive negator and {$H_{\eta}=F^{\downarrow,{\eta}}_{\mathcal{P}}$}.Then for all $y\in Y,x\in X$
{\begin{eqnarray*}
({\mathbf{N}}(H_{\eta}[{\mathbf{N}}( 1_{\lbrace x\rbrace})]))(y)&=& ({\mathbf{N}}( F^{\downarrow,{\eta}}_{\mathcal{P}}[{\mathbf{N}}( 1_{\lbrace x\rbrace})]))(y)\\
&=&{\mathbf{N}}( F^{\downarrow,{\eta}}_{\mathcal{P}}[{\mathbf{N}}( 1_{\lbrace x\rbrace})](y))\\
&=& {\mathbf{N}}(\bigwedge_{z\in X}{\eta}(\mathbf{N}_{{\eta}}(A_{y}(z)),({\mathbf{N}}( 1_{\lbrace x\rbrace}))(z)))\\
&=& {\mathbf{N}}(\bigwedge_{z\in X}{\eta}(\mathbf{N}(A_{y}(z)),{\mathbf{N}}( 1_{\lbrace x\rbrace}(z))))\\
&=&{\mathbf{N}}({\eta}(\mathbf{N}(A_{y}(x)),0))\\
&=&{\mathbf{N}}(\mathbf{N}(A_{y}(x)))\\
&=&A_{y}(x).
\end{eqnarray*}
 Thus $ ({\mathbf{N}}(H_{\eta}[{\mathbf{N}}( 1_{\lbrace x\rbrace})]))(y)=1$ iff $A_{y}(x)=1$ iff $v(x)=y$.} From Propositions \ref{p35} and \ref{p36}, $(X,Y,v,H_{\eta})$ is an {$L$-fuzzy} lower transformation system on $X$ {determined by ${\eta}$}.
\begin{pro}\label{p54}
Let $\mathcal{I}_{{\theta}}$ be an $EP$-residual implicator over $L$ such that $\mathbf{N}_{\mathcal{I}_{{\theta}}}$ is an involutive negator. Then the following statements are equivalent:
\begin{itemize}
\item[(i)] $\mathcal{H}_{\mathcal{I}_{{\theta}}}=(X,Y, v, H_{\mathcal{I}_{{\theta}}})$ is an $L$-fuzzy lower transformation system on $X$ determined by $\mathcal{I}_{{\theta}}$ and {$Y\subseteq X$}.
\item[(ii)] There exists an $L$-fuzzy partition $\mathcal{P}$ of $X$ {indexed by $Y$}, such that $v(x) = y$ iff $x\in core(A_y)$ and $H_{\mathcal{I}_{{\theta}}} = F^{\downarrow,\mathcal{I}_{{\theta}}}_{\mathcal{P}}$.
\end{itemize}
\end{pro}
\textbf{Proof:} Let ${\mathcal{H}_{\mathcal{I}_{{\theta}}}=(X,Y,v,H_{\mathcal{I}_{{\theta}}})}$ be an {$L$-fuzzy} lower transformation system on X {determined by $\mathcal{I}_{{\theta}}$}. Also, let $\mathcal{P} = \lbrace A_y : y\in Y\rbrace$ such that for all $y\in Y$, $A_y \in L^X$ is given by {$A_y(x) = ({\mathbf{N}_{\mathcal{I}_{{\theta}}}}( H_{\mathcal{I}_{{\theta}}}[{\mathbf{N}_{\mathcal{I}_{{\theta}}}}( 1_{\lbrace x\rbrace})]))(y)$}, $x\in X$. Now, from Definition \ref{UTS}(iii), $A_{v(x)}(x) = {({\mathbf{N}_{\mathcal{I}_{{\theta}}}}( H_{\mathcal{I}_{{\theta}}}[{\mathbf{N}_{\mathcal{I}_{{\theta}}}}( 1_{\lbrace x\rbrace})]))(v(x))} = 1$, or that, $x \in core( A_{v(x)})$. Further, {for $t \in core ( A_y) \cap core (A_z)$, $y, z \in Y$ and the fact that ${\mathbf{N}_{\mathcal{I}_{{\theta}}}}(x)=\mathcal{I}_{\theta}(x,0)$}, ${({\mathbf{N}_{\mathcal{I}_{{\theta}}}}( H_{\mathcal{I}_{{\theta}}}[{\mathbf{N}_{\mathcal{I}_{{\theta}}}}( 1_{\lbrace t\rbrace})]))(y)= 1 =({\mathbf{N}_{\mathcal{I}_{{\theta}}}}( H_{\mathcal{I}_{{\theta}}}[{\mathbf{N}_{\mathcal{I}_{{\theta}}}}( 1_{\lbrace t\rbrace})]))(z)}$, i.e., $ A_y(t) = 1 = A_z(t)$ iff $y = v(t) = z$. Thus $\lbrace core(A_y):y\in Y\rbrace$ is a partition of $X$ and therefore $\mathcal{P}$ is an {$L$-fuzzy} partition of $X$. Now, for all $y\in Y$ and $f\in L^X$
\begin{eqnarray*}
F^{\downarrow,\mathcal{I}_{{\theta}}} _{\mathcal{P}}[f](y)&=&\bigwedge\limits_{x\in X}\mathcal{I}_{\theta}( A_y(x),f(x))\\
&=&\bigwedge\limits_{x\in X}\mathcal{I}_{\theta}( {({\mathbf{N}_{\mathcal{I}}}(H_{\mathcal{I}_{{\theta}}}[{\mathbf{N}_{\mathcal{I}_{{\theta}}}}( 1_{\lbrace x\rbrace})]))(y)}, f(x))\\
&=& {{\bigwedge\limits_{x\in X}\mathcal{I}_{\theta}({\mathbf{N}_{\mathcal{I}_{{\theta}}}}}(H_{\mathcal{I}_{{\theta}}}[{\mathbf{N}_{\mathcal{I}_{{\theta}}}}( 1_{\lbrace x\rbrace})](y)), {\mathbf{N}_{\mathcal{I}_{{\theta}}}}({\mathbf{N}_{\mathcal{I}_{{\theta}}}}(f(x))))}\\
&=&\bigwedge\limits_{x\in X}\mathcal{I}_{\theta}({{\mathbf{N}_{\mathcal{I}_{{\theta}}}}}( f(x)), {{\mathbf{N}_{\mathcal{I}_{{\theta}}}}}({\mathbf{N}_{\mathcal{I}_{{\theta}}}}( H_{\mathcal{I}_{{\theta}}}[{{\mathbf{N}_{\mathcal{I}_{{\theta}}}}}(1_{\lbrace x\rbrace})](y))))\\
&=&{\bigwedge\limits_{x\in X}\mathcal{I}_{\theta}(({\mathbf{N}_{\mathcal{I}_{{\theta}}}}(f))(x), H_{\mathcal{I}_{{\theta}}}[{\mathbf{N}_{\mathcal{I}_{{\theta}}}}( 1_{\lbrace x\rbrace})](y))}\\
&=&{\bigwedge\limits_{x\in X}H_{\mathcal{I}_{{\theta}}}[{\mathcal{I}_{{\theta}}}({{\mathbf{N}_{\mathcal{I}_{{\theta}}}}(\textbf{f(x)})}, {\mathbf{N}_{\mathcal{I}_{{\theta}}}}( 1_{\lbrace x\rbrace}))](y)}\\
&=&{H_{\mathcal{I}_{{\theta}}}[\bigwedge\limits_{x\in X}{\mathcal{I}_{{\theta}}}({{\mathbf{N}_{\mathcal{I}_{{\theta}}}}(\textbf{f(x)})}, {\mathbf{N}_{\mathcal{I}_{{\theta}}}}( 1_{\lbrace x\rbrace}))](y)}\\
&=&H_{\mathcal{I}_{{\theta}}}[{\mathbf{N}_{\mathcal{I}_{{\theta}}}}({\mathbf{N}_{\mathcal{I}_{{\theta}}}}( f))](y)\\
&=&H_{\mathcal{I}_{{\theta}}}[f](y).
\end{eqnarray*}
Thus $H_{\mathcal{I}_{{\theta}}}=F^{\downarrow,\mathcal{I}_{{\theta}}}_{\mathcal{P}}$.
Conversely, let $\mathcal{P}=\{A_y\in L^X: y\in Y\}$ be an {$L$-fuzzy} partition of base set $X\neq\emptyset$. Let us define a map $v: X \rightarrow Y$ such that $v(x) = y$ iff $x \in core(A_{y})$. Further, let $\mathcal{I}_{{\theta}}$ be a residual implicator such that ${\mathbf{N}_{\mathcal{I}_{{\theta}}}}(\cdot)=\mathcal{I}_{\theta}(\cdot,0)$ is an involutive negator) and {$H_{\mathcal{I}_{{\theta}}}=F^{\downarrow,\mathcal{I}_{{\theta}}}_{\mathcal{P}}$}. Then for all $y\in Y,x\in X$
{\begin{eqnarray*}
({\mathbf{N}_{\mathcal{I}_{{\theta}}}}(H_{\mathcal{I}_{{\theta}}}[{\mathbf{N}_{\mathcal{I}_{{\theta}}}}( 1_{\lbrace x\rbrace})]))(y)&=& ({\mathbf{N}_{\mathcal{I}_{{\theta}}}}( F^{\downarrow,\mathcal{I}_{{\theta}}}_{\mathcal{P}}[{\mathbf{N}_{\mathcal{I}_{{\theta}}}}( 1_{\lbrace x\rbrace})]))(y)\\
&=&{\mathbf{N}_{\mathcal{I}_{{\theta}}}}( F^{\downarrow,\mathcal{I}_{{\theta}}}_{\mathcal{P}}[{\mathbf{N}_{\mathcal{I}_{{\theta}}}}( 1_{\lbrace x\rbrace})](y))\\
&=& {\mathbf{N}_{\mathcal{I}_{{\theta}}}}(\bigwedge_{z\in X}\mathcal{I}_{\theta}(A_{y}(z),({\mathbf{N}_{\mathcal{I}_{{\theta}}}}( 1_{\lbrace x\rbrace}))(z)))\\
&=& {\mathbf{N}_{\mathcal{I}_{{\theta}}}}(\bigwedge_{z\in X}\mathcal{I}_{\theta}(A_{y}(z),{\mathbf{N}_{\mathcal{I}_{{\theta}}}}( 1_{\lbrace x\rbrace}(z))))\\
&=&{\mathbf{N}_{\mathcal{I}_{{\theta}}}}(\mathcal{I}_{\theta}(A_{y}(x),0))=A_{y}(x).
\end{eqnarray*}
 Thus $ ({\mathbf{N}_{\mathcal{I}_{{\theta}}}}(H_{\mathcal{I}_{{\theta}}}[{\mathbf{N}_{\mathcal{I}_{{\theta}}}}( 1_{\lbrace x\rbrace})]))(y)=1$ iff $A_{y}(x)=1$ iff $v(x)=y$.} From Propositions \ref{p35} and \ref{p36}, $(X,Y,v,H_{\mathcal{I}_{{\theta}}})$ is an {$L$-fuzzy} lower transformation system on $X$ {determined by ${\mathcal{I}_{{\theta}}}$}.\\\\
 Next, we have the following.
\begin{pro}\label{p55}
Let ${\theta}$ and ${\eta}$ be dual with respect to an involutive negator $\mathbf{N}$, $\mathcal{U}_{{\theta}}=(X,Y,u,U_{{\theta}})$ and $\mathcal{H}_{\eta}=(X,Y,u,H_{\eta})$ be $L$-fuzzy upper and lower transformation systems, respectively. Then there exists an $L$-fuzzy partition $\mathcal{P}$ such that $U_{{\theta}}=F^{\uparrow,{\theta}}_{\mathcal{P}},\linebreak H_{\eta}=F_{\mathcal{P}}^{\downarrow,{\eta}}$ iff for all $f\in L^X$,
\begin{itemize}
\item[(i)] $U_{{\theta}}[f]=\mathbf{N}(H_{{\eta}}[{\mathbf{N}}(f)])$, i.e., $\mathbf{N}(U_{{\theta}}[f])=H_{\eta}[{\mathbf{N}}(f)]$, and
\item[(ii)] $H_{\eta}[f]=\mathbf{N}(U_{{\theta}}[{\mathbf{N}}(f)])$, i.e, $\mathbf{N}(H_{\eta})=U_{{\theta}}[{\mathbf{N}}(f)]$.
\end{itemize}
\end{pro}
\textbf{Proof:} From Propositions \ref{p32a}, is can be easily show that conditions (i) and (ii) hold. Now, we only need to show that the converse part. For which, let condition (i) holds. Further, let $\{A_{1,y}:y\in Y\},\{A_{2,y}:y\in Y\}\subseteq L^X$ such that $A_{1,y}(x)=U_{{\theta}}[1_{\{x\}}](y),\linebreak A_{2,y}(x)=\mathbf{N}(H_{\eta}[\mathbf{N}(1_{\{x\}})])(y),\,\forall\,x\in X,y\in Y$. Then from propositions \ref{p51} and \ref{p52}, it is clear that $\{A_{1,y}:y\in Y\},\{A_{2,y}:y\in Y\}\subseteq L^X$ are $L$-fuzzy partitions of $X$ and $U_{{\theta}}=F^{\uparrow,{\theta}}_{1,\mathcal{P}},H_{\eta}=F_{2,\mathcal{P}}^{\downarrow,{\eta}}$. Now, from condition (i), we have $U_{{\theta}}[f]=\mathbf{N}(H_{{\eta}}[{\mathbf{N}}(f)])=\mathbf{N}(F_{2,\mathcal{P}}^{\downarrow,{\eta}}[{\mathbf{N}}(f)])=F^{\uparrow,{\theta}}_{2,\mathcal{P}}[f]$. Thus $F^{\uparrow,{\theta}}_{1,\mathcal{P}}=F^{\uparrow,{\theta}}_{2,\mathcal{P}}$ and $A_{1y}=A_{2y},\,\forall\,y\in Y$. Similarly, we can show that when condition (ii) holds.
 \begin{pro}\label{p46}
Let ${\theta}$ and ${\eta}$ be dual with respect to an involutive negator $\mathbf{N}$, $\mathcal{U}_{\mathcal{I}_{\eta}}=(X,Y,u,U_{\mathcal{I}_{\eta}})$ and $\mathcal{H}_{\mathcal{I}_{{\theta}}}=(X,Y,u,H_{\mathcal{I}_{\eta}})$ be $L$-fuzzy upper and lower transformation systems, respectively. Then there exists an $L$-fuzzy partition $\mathcal{P}$ such that $F^{\uparrow,_{\mathcal{I}_{\eta}}}_{\mathcal{P}}=U_{\mathcal{I}_{\eta}},F_{\mathcal{P}}^{\downarrow,_{\mathcal{I}_{{\theta}}}}=H_{\mathcal{I}_{{\theta}}}$ iff for all $f\in L^X$
\begin{itemize}
\item[(i)] $U_{\mathcal{I}_{\eta}}[f]=\mathbf{N}(H_{\mathcal{I}_{{\theta}}}[{\mathbf{N}}(f)])$, i.e., $\mathbf{N}(U_{\mathcal{I}_{\eta}}[f])=H_{\mathcal{I}_{{\theta}}}[{\mathbf{N}}(f)]$, and
\item[(ii)] $H_{\mathcal{I}_{{\theta}}}[f]=\mathbf{N}(U_{\mathcal{I}_{\eta}}[{\mathbf{N}}(f)])$, i.e, $\mathbf{N}(H_{\mathcal{I}_{{\theta}}})=U_{\mathcal{I}_{\eta}}[{\mathbf{N}}(f)]$.
\end{itemize}
\end{pro}
\textbf{Proof:} Similar to that of Proposition \ref{p46}.
 \begin{pro}\label{p47}
Let $\mathbf{N}_{\mathcal{I}_{{\theta}}}$ be involutive negator, $\mathcal{U}_{{\theta}}=(X,Y,u,U_{{\theta}})$ and $\mathcal{H}_{\mathcal{I}_{{\theta}}}=\linebreak(X,Y,u,H_{\mathcal{I}_{\eta}})$ be $L$-fuzzy upper and lower transformation systems, respectively. Then there exists an $L$-fuzzy partition $\mathcal{P}$ such that $F^{\uparrow,_{{\theta}}}_{\mathcal{P}}=U_{{\theta}},F_{\mathcal{P}}^{\downarrow,_{\mathcal{I}_{{\theta}}}}=H_{\mathcal{I}_{{\theta}}}$ iff 
\begin{itemize}
\item[(i)] $U_{{\theta}}[f]=\mathbf{N}_{\mathcal{I}_{{\theta}}}(H_{\mathcal{I}_{{\theta}}}[{\mathbf{N}_{\mathcal{I}_{{\theta}}}}(f)])$, i.e., $\mathbf{N}_{\mathcal{I}_{{\theta}}}(U_{{\theta}}[f])=H_{\mathcal{I}_{{\theta}}}[{\mathbf{N}_{\mathcal{I}_{{\theta}}}}(f)]$, and
\item[(ii)] $H_{\mathcal{I}_{{\theta}}}[f]=\mathbf{N}_{\mathcal{I}_{{\theta}}}(U_{{\theta}}[{\mathbf{N}_{\mathcal{I}_{{\theta}}}}(f)])$, i.e, $\mathbf{N}_{\mathcal{I}_{{\theta}}}(H_{\mathcal{I}_{{\theta}}})=U_{{\theta}}[{\mathbf{N}_{\mathcal{I}_{{\theta}}}}(f)]$.
\end{itemize}
\end{pro}
\textbf{Proof:} Similar to that of Proposition \ref{p46}.
 \begin{pro}\label{p48}
Let $\mathbf{N}_{\mathcal{I}_{\eta}}$ be involutive negators, $\mathcal{U}_{\mathcal{I}_{\eta}}=(X,Y,u,U_{\mathcal{I}_{\eta}})$ and $\mathcal{H}_{{\eta}}=(X,Y,u,H_{{\eta}})$ be $L$-fuzzy upper and lower transformation systems, respectively. Then there exists an $L$-fuzzy partition $\mathcal{P}$ such that $F^{\uparrow,\mathcal{I}_{\eta}}_{\mathcal{P}}=U_{\mathcal{I}_{\eta}},F_{\mathcal{P}}^{\downarrow,{\eta}}=H_{{\eta}}$ iff for all $f\in L^X$
\begin{itemize}
\item[(i)] $U_{\mathcal{I}_{\eta}}[f]=\mathbf{N}_{\mathcal{I}_{\eta}}(H_{{\eta}}[{\mathbf{N}_{\mathcal{I}_{\eta}}}(f)])$, i.e., $\mathbf{N}_{\mathcal{I}_{\eta}}(U_{\mathcal{I}_{\eta}}[f])=H_{{\eta}}[{\mathbf{N}_{\mathcal{I}_{\eta}}}(f)]$, and
\item[(ii)] $H_{{\eta}}[f]=\mathbf{N}_{\mathcal{I}_{\eta}}(U_{\mathcal{I}_{\eta}}[{\mathbf{N}_{\mathcal{I}_{\eta}}}(f)])$, i.e, $\mathbf{N}_{\mathcal{I}_{\eta}}(H_{{\eta}})=U_{\mathcal{I}_{\eta}}[{\mathbf{N}_{\mathcal{I}_{\eta}}}(f)]$.
\end{itemize}
\end{pro}
\textbf{Proof:} Similar to that of Proposition \ref{p46}.
\section{Concluding remarks}
In this contribution, we have presented the theory of direct $F$-transforms determined by overlap and grouping maps, residual and co-residual implicators from both constructive and axiomatic approaches. In which, $F^{\uparrow,{\theta}},F^{\downarrow,{\eta}},F^{\downarrow,\mathcal{I}_{{\theta}}}$ are the extension of the direct $F$-transforms introduced in \cite{per,anan,tri} and $F^{\uparrow,\mathcal{I}_{\eta}}$ is a new definition. The main contributions of this paper are listed as follows.
\begin{itemize}
\item[$\bullet$] We have shown the duality of the proposed direct $F$-transform and established a connection among these direct $F$-transforms. In addition, we have discussed the basic results of these direct $F$-transforms.
\item[$\bullet$] We have introduced the idea of the inverse of these $F$-transforms. Further, we have shown that the original $L$-fuzzy set and inverse of these $F$-transform have the same $F$-transform under certain conditions. 
\item[$\bullet$] Further, we have shown an axiomatic characterization of the proposed direct $F$-transforms. 
\item[$\bullet$] Finally, the duality of $L$-fuzzy transformation systems has been examined.
\end{itemize}
Both the theories viz., theory of $F$-transforms and the theory of overlap and grouping maps have already shown to be helpful in practical applications. Accordingly, combining both ideas may provide us with new applications in data analysis and image processing problems.

\end{document}